\documentstyle[11pt,amstex]{article}

\font\tenfrak=eufm10
\font\sevenfrak=eufm7
\font\fivefrak=eufm5
\newfam\frakfam \def\frak{\fam\frakfam\tenfrak} \textfont\frakfam=\tenfrak
  \scriptfont\frakfam=\sevenfrak  \scriptscriptfont\frakfam=\fivefrak

\begin{document}
    \pagestyle{plain}
    \setlength{\baselineskip}{1.2\baselineskip}
    \setlength{\parindent}{\parindent}

\title{{\bf Classical dynamical $r$-matrices and homogeneous Poisson 
structures on $G/H$ and $K/T$}}
\author{Jiang-Hua Lu
\thanks{Research partially supported by
an NSF Postdoctorial Fellowship and by NSF grant DMS 9803624.}\\
Department of Mathematics, University of Arizona, Tucson, AZ 85721\\
}
\maketitle

\newtheorem{thm}{Theorem}[section]
\newtheorem{lem}[thm]{Lemma}
\newtheorem{prop}[thm]{Proposition}
\newtheorem{cor}[thm]{Corollary}
\newtheorem{rem}[thm]{Remark}
\newtheorem{exam}[thm]{Example}
\newtheorem{nota}[thm]{Notation}
\newtheorem{dfn}[thm]{Definition}
\newtheorem{ques}[thm]{Question}
\newtheorem{eq}{thm}

\newcommand{\rw}{\rightarrow}
\newcommand{\lrw}{\longrightarrow}
\newcommand{\rhu}{\rightharpoonup}
\newcommand{\lhu}{\leftharpoonup}
\newcommand{\Map}{\longmapsto}
\newcommand{\qed}{\begin{flushright} {\bf Q.E.D.}\ \ \ \ \
                  \end{flushright} }
\newcommand{\beqa}{\begin{eqnarray*}}
\newcommand{\eeqa}{\end{eqnarray*}}

\newcommand{\la}{\mbox{$\langle$}}
\newcommand{\ra}{\mbox{$\rangle$}}
\newcommand{\ot}{\mbox{$\otimes$}}
\newcommand{\xa}{\mbox{$x_{(1)}$}}
\newcommand{\xb}{\mbox{$x_{(2)}$}}
\newcommand{\xc}{\mbox{$x_{(3)}$}}
\newcommand{\ya}{\mbox{$y_{(1)}$}}
\newcommand{\yb}{\mbox{$y_{(2)}$}}
\newcommand{\yc}{\mbox{$y_{(3)}$}}
\newcommand{\yd}{\mbox{$y_{(4)}$}}
\renewcommand{\aa}{\mbox{$a_{(1)}$}}
\newcommand{\ab}{\mbox{$a_{(2)}$}}
\newcommand{\ac}{\mbox{$a_{(3)}$}}
\newcommand{\ad}{\mbox{$a_{(4)}$}}
\newcommand{\ba}{\mbox{$b_{(1)}$}}
\newcommand{\bt}{\mbox{$b_{(2)}$}}
\newcommand{\bc}{\mbox{$b_{(3)}$}}
\newcommand{\ca}{\mbox{$c_{(1)}$}}
\newcommand{\cb}{\mbox{$c_{(2)}$}}
\newcommand{\cc}{\mbox{$c_{(3)}$}}
\newcommand{\uo}{\mbox{$\underbar{o}$}}

\newcommand{\ts}{\mbox{$\sigma$}}
\newcommand{\las}{\mbox{${}_{\sigma}\!A$}}
\newcommand{\lasone}{\mbox{${}_{\sigma'}\!A$}}
\newcommand{\ras}{\mbox{$A_{\sigma}$}}
\newcommand{\rds}{\mbox{$\cdot_{\sigma}$}}
\newcommand{\lds}{\mbox{${}_{\sigma}\!\cdot$}}

\newcommand{\bb}{\mbox{$\bar{\beta}$}}
\newcommand{\bg}{\mbox{$\bar{\gamma}$}}

\newcommand{\id}{\mbox{${\em id}$}}
\newcommand{\Fun}{\mbox{${\em Fun}$}}
\newcommand{\End}{\mbox{${\em End}$}}
\renewcommand{\ker}{\mbox{${\rm Ker}$}}
\newcommand{\im}{\mbox{${\rm Im}$}}

\newcommand{\Hom}{\mbox{${\em Hom}$}}
\newcommand{\ta}{\mbox{${\mbox{$\scriptscriptstyle A$}}$}}
\newcommand{\ms}{\mbox{${\mbox{$\scriptscriptstyle M$}}$}}
\newcommand{\ap}{\mbox{$A_{\mbox{$\scriptscriptstyle P$}}$}}
\newcommand{\tx}{\mbox{$\mbox{$\scriptscriptstyle X$}$}}
\newcommand{\txo}{\mbox{$\mbox{$\scriptscriptstyle X_1$}$}}
\newcommand{\tE}{\mbox{$\mbox{$\scriptscriptstyle E$}$}}
\newcommand{\ty}{\mbox{$\mbox{$\scriptscriptstyle Y$}$}}
\newcommand{\kt}{\mbox{$K_{\tx}$}}
\newcommand{\pk}{\mbox{$\pi_{\scriptscriptstyle X_1,
          \lambda}^{\scriptscriptstyle X}$}}
\newcommand{\kk}{\mbox{$K \times_{\scriptscriptstyle K_{X}} (K_{\tx}/T)$}}
\newcommand{\axl}{\mbox{$A_{\tx}(\lambda)$}}
\newcommand{\pxlp}{\mbox{$\partial^{'}_{\tx, \lambda}$}}
\newcommand{\pxlpp}{\mbox{$\partial^{''}_{\tx, \lambda}$}}
\newcommand{\sxl}{\mbox{$S_{\tx, \lambda}$}}
\newcommand{\tpxlp}{\mbox{$\tilde{\partial}^{'}_{\tx, \lambda}$}}
\newcommand{\tpxlpp}{\mbox{$\tilde{\partial}^{''}_{\tx, \lambda}$}}
\newcommand{\tpxl}{\mbox{$\tilde{\partial}_{\tx, \lambda}$}}
\newcommand{\gxl}{\mbox{$g_{\tx, \lambda}$}}
\newcommand{\hxl}{\mbox{$h_{\tx, \lambda}$}}
\newcommand{\lxl}{\mbox{$L_{\tx, \lambda}$}}
\newcommand{\tspxl}{\mbox{$\tilde{\partial}_{\ast, \tx, \lambda}$}}
\newcommand{\ixl}{\mbox{$I_{\tx, \lambda}$}}
\newcommand{\mxl}{\mbox{$m_{\tx, \lambda}$}}
\newcommand{\tmxl}{\mbox{$\tilde{m}_{\tx, \lambda}$}}
\newcommand{\starxl}{\mbox{$\star_{\tx, \lambda}$}}
\newcommand{\astxl}{\mbox{$\ast_{\tx, \lambda}$}}
\newcommand{\bxl}{\mbox{$b_{\tx, \lambda}$}}

\newcommand{\wm}{W^{\tx}}
\newcommand{\wx}{W_{\tx}}
\newcommand{\dw}{\dot{w}}
\newcommand{\sw}{\Sigma_w}
\newcommand{\swa}{\Sigma_{w_1}}
\newcommand{\swb}{\Sigma_{w_1 w_2}}
\newcommand{\cw}{C_{\dot{w}}}
\newcommand{\cwa}{C_{\dot{w_1}}}

\newcommand{\pp}{\mbox{$\pi_{\mbox{$\scriptscriptstyle P$}}$}}
\newcommand{\pg}{\mbox{$\pi_{\mbox{$\scriptscriptstyle G$}}$}}
\newcommand{\pge}{\mbox{$\pi_{\Sigma_{+}}$}}

\newcommand{\Xa}{\mbox{$X_{\alpha}$}}
\newcommand{\Ya}{\mbox{$Y_{\alpha}$}}
\newcommand{\pix}{\mbox{$\pi_{\tx, \txo, \lambda}$}}
\newcommand{\qix}{\mbox{$\pi_{\tx, \emptyset, \lambda}$}}

\newcommand{\lix}{\mbox{$\fl_{\tx, \txo, \lambda}$}}
\newcommand{\liy}{\mbox{$\fl_{\ty, \txo, \lambda}$}}
\newcommand{\lio}{\mbox{$\fl_{\txo, \lambda}$}}

\newcommand{\pinf}{\mbox{$\pi_{\infty}$}}
\newcommand{\fix}{\mbox{$\fl_{\tx,\lambda}$}}
\newcommand{\bix}{\mbox{$b_{\tx,\lambda}$}}
\newcommand{\dix}{\mbox{$d_{\tx,\lambda}$}}
\newcommand{\pxl}{\mbox{$\partial_{\tx, \lambda}$}}
\newcommand{\tml}{\mbox{$\tilde{M}_{\lambda}$}}
\newcommand{\tse}{\mbox{$T^{*}_{e}(K/T)$}}
\newcommand{\te}{\mbox{$T_{e}(K/T)$}}

\newcommand{\el}{\mbox{$e^{- \lambda}$}}
\newcommand{\pkl}{\mbox{$\pi^{{\scriptscriptstyle K_{X}}}_{\lambda}$}}

\newcommand{\dpi}{\mbox{$\delta_{\pi}$}}

\newcommand{\asemi}{\mbox{$\ap \#_{\sigma} A^*$}}
\newcommand{\dsemi}{\mbox{$A \#_{\Delta} A^*$}}

\newcommand{\semi}{\mbox{$\times_{{\frac{1}{2}}}$}}
\newcommand{\fd}{\mbox{${\frak d}$}}
\newcommand{\fa}{\mbox{${\frak a}$}}
\newcommand{\ft}{\mbox{${\frak t}$}}
\newcommand{\fk}{\mbox{${\frak k}$}}
\newcommand{\fg}{\mbox{${\frak g}$}}
\newcommand{\fq}{\mbox{${\frak q}$}}
\newcommand{\fl}{\mbox{${\frak l}$}}
\newcommand{\fs}{\mbox{${\frak s}$}}

\newcommand{\flp}{\mbox{${\frak l}_{p}$}}
\newcommand{\fh}{\mbox{${\frak h}$}}
\newcommand{\fn}{\mbox{${\frak n}$}}
\newcommand{\fm}{\mbox{${\frak m}$}}
\newcommand{\fu}{{\frak u}}

\newcommand{\fks}{\mbox{$\fk_{\sigma}$}}
\newcommand{\Ks}{\mbox{$K_{\sigma}$}}
\newcommand{\Ps}{\mbox{$\Pi_{\sigma}$}}
\newcommand{\ps}{\mbox{$\pi_{\sigma}$}}

\newcommand{\bfny}{\mbox{${\bar{\frak n}}_{\ty}$}}

\newcommand{\fp}{\mbox{${\frak p}$}}
\newcommand{\fb}{\mbox{${\frak b}$}}
\newcommand{\fbp}{\mbox{${\frak b}_{+}$}}
\newcommand{\fbm}{\mbox{${\frak b}_{-}$}}
\newcommand{\fnp}{\mbox{${\frak n}_{+}$}}
\newcommand{\fnm}{\mbox{${\frak n}_{-}$}}
\newcommand{\fgs}{\mbox{${\frak g}^*$}}
\newcommand{\wg}{\mbox{$\wedge {\frak g}$}}
\newcommand{\wgs}{\mbox{$\wedge {\frak g}^*$}}
\newcommand{\wxl}{\mbox{$x_1 \wedge x_2 \wedge \cdots \wedge x_l$}}
\newcommand{\wxk}{\mbox{$x_1 \wedge x_2 \wedge \cdots \wedge x_k$}}
\newcommand{\wyl}{\mbox{$y_1 \wedge y_2 \wedge \cdots \wedge y_l$}}
\newcommand{\wxkm}{\mbox{$x_1 \wedge x_2 \wedge \cdots \wedge x_{k-1}$}}
\newcommand{\wxik}{\mbox{$\xi_1 \wedge \xi_2 \wedge \cdots \wedge \xi_k$}}
\newcommand{\wxikm}{\mbox{$\xi_1 \wedge \cdots \wedge \xi_{k-1}$}}
\newcommand{\wetal}{\mbox{$\eta_1 \wedge \eta_2 \wedge \cdots \wedge \eta_l$}}

\newcommand{\winv}{\mbox{$(\wedge \fg_{1}^{\perp})^{\fg_1}$}}
\newcommand{\wetak}{\mbox{$\eta_1 \wedge \cdots \wedge \eta_k$}}
\newcommand{\gonep}{\mbox{$\fg_{1}^{\perp}$}}
\newcommand{\wonep}{\mbox{$\wedge \fg_{1}^{\perp}$}}

\newcommand{\db}{\mbox{$\fd = \fg \bowtie \fgs$}}
\newcommand{\fds}{\mbox{${\scriptscriptstyle {\frak d}}$}}

\newcommand{\Gs}{\mbox{$G^*$}}
\newcommand{\pis}{\mbox{$\pi_{\sigma}$}}
\newcommand{\ea}{\mbox{$E_{\alpha}$}}
\newcommand{\eb}{\mbox{$E_{-\alpha}$}}
\newcommand{\Bm}{\mbox{$ {}^B \! M$}}
\newcommand{\kBm}{\mbox{$ {}^B \! M^k$}}
\newcommand{\Bb}{\mbox{$ {}^B \! b$}}
\newcommand{\epe}{\mbox{$\varepsilon$}}
\newcommand{\eot}{\mbox{${\epe \over 2}$}}

\newcommand{\cfg}{\mbox{$C(\fg \oplus \fgs)$}}
\newcommand{\backl}{\mathbin{\vrule width1.5ex height.4pt\vrule height1.5ex}}
 
\newcommand{\bx}{\mbox{${\bar{x}}$}}
\newcommand{\by}{\mbox{${\bar{y}}$}}
\newcommand{\bz}{\mbox{${\bar{z}}$}}
\newcommand{\pgs}{\mbox{${\pi_{\mbox{\tiny G}^{*}}}$}}

\newcommand{\tk}{\mbox{${\scriptscriptstyle K}$}}
\newcommand{\tg}{\mbox{${\scriptscriptstyle G}$}}

\newcommand{\piK}{\mbox{$\pi_{\tk}$}}
\newcommand{\piG}{\mbox{$\pi_{\tg}$}}

\newcommand{\piKX}{\mbox{$\pi_{\scriptscriptstyle K_X}$}}

\newcommand{\tlp}{\mbox{$\tilde{\pi}$}}
\newcommand{\tlpz}{\mbox{$\tilde{\pi}_0$}}

\newcommand{\hlp}{\mbox{$\hat{\pi}$}}
\newcommand{\tp}{\mbox{$\tilde{\pi}$}}
\newcommand{\sn}{\mbox{$s_{\scriptscriptstyle N}$}}
\newcommand{\tn}{\mbox{$t_{\scriptscriptstyle N}$}}
\newcommand{\sm}{\mbox{$s_{\scriptscriptstyle M}$}}
\newcommand{\tm}{\mbox{$t_{\scriptscriptstyle M}$}}

\newcommand{\C}{\mbox{${\Bbb C}$}}
\newcommand{\Z}{\mbox{${\Bbb Z}$}}

\newcommand{\R}{\mbox{${\Bbb R}$}}
\renewcommand{\a}{\mbox{$\alpha$}}

\newcommand{\Lagr}{{\cal L}}
\newcommand{\su}{{\frak s}{\frak u}}
\newcommand{\fsl}{{\frak s}{\frak l}}
\newcommand{\Gr}{{\rm Gr}}

\vspace{-.2in}
\begin{abstract}
Let $G$ be a finite dimensional simple complex group 
equipped with the 
standard Poisson Lie group structure. We show that all $G$-homogeneous 
(holomorphic) Poisson structures on $G/H$, where $H \subset G$
is a Cartan subgroup, come from solutions 
to the Classical Dynamical Yang-Baxter equations which are classified 
by Etingof and Varchenko. A similar result holds for 
the maximal compact subgroup $K$, and we get a family 
of $K$-homogeneous Poisson structures on $K/T$, where $T = K \cap H$
is a maximal torus of $K$. 
This family
exhausts all  $K$-homogeneous Poisson structures on $K/T$ up to
isomorphisms. We study some Poisson geometrical
properties of members of this family such as their
symplectic leaves, their modular classes, and the moment maps
for the $T$-action.
\end{abstract}

\tableofcontents

\section{Introduction}
\label{sec_intro}

This paper is motivated by the work of Etingof
and Varchenko \cite{e-v:cdyb} on 
{\it classical  dynamical $r$-matrices} for the
pair $(\fg, \fh)$, where $\fg$ is a complex simple Lie algebra
and $\fh \subset \fg$ a Cartan subalgebra. 

A classical dynamical $r$-matrix is, by definition, a 
meromorphic function $r: \fh^* \rightarrow \fg \ot \fg$
satisfying the so-called {\it
Classical Dynamical Yang-Baxter Equation} (CDYBE):
\[
{\rm Alt} (dr) \, + \, 
[r^{12},r^{13}] + [r^{12}, r^{23}] + [r^{13}, r^{23}] \, = \, 0.
\]
(See Section \ref{sec_cdybe} for details). One such $r$-matrix
has the form
\[
r(\lambda) \, = \, {\epe \over 2}\, \Omega \,+ \, 
{\epe \over 2} \sum _{\alpha \in \Sigma} \,
\coth ({\epe \over 2} \ll \alpha, \lambda \gg) \ea \otimes \eb,
\]
where $\Omega \in (S^2 \fg)^{\frak g}$ corresponds to 
the Killing form $\ll \, , \, \gg$ of $\fg, \Sigma$ is
the set of roots of $\fg$ with respect to $\fh$, the
$\ea$ and $\eb$'s are root vectors, and 
$\coth (x)  =  {e^x + e^{-x} \over e^x - e^{-x}}$
is the hyperbolic cotangent function. Other $r$-matrices
can be obtained by performing certain ``gauge transformations" to 
the one above and by taking various limits of it.
See Section \ref{sec_cdybe}.

We wanted to understand the geometrical meaning of these
$r$-matrices. Etingof and Varchenko show in \cite{e-v:cdyb}
that every classical dynamical $r$-matrix defines a Poisson
groupoid over an open subset of $\fh^*$.  In this paper, we
give another geometrical interpretation of the $r$-matrices
by connecting them with Poisson structures on the spaces $G/H$ 
and $K/T$, where
$G$ is a complex Lie group with Lie algebra $\fg$, $H \subset G$
its connected subgroup corresponding to $\fh$, $K$ a compact real
form of $G$, and $T = K \cap H$. We then study some Poisson
geometrical properties of these Poisson structures on $K/T$ such as
their symplectic leaves, their modular classes, and the moment
maps for the $T$-action. 
We now explain this in more detail.

A special example of a classical dynamical $r$-matrix is one that 
is not ``dynamical", i.e., independent 
of $\lambda$. It is given by
\[
r_0 \, = \, {\frac{\epe}{2}} \Omega \, + \, 
 c \, + \, {\frac{\epe}{2}} \sum_{\alpha \in \Sigma_{+}}
E_{\alpha} \wedge E_{-\alpha}
\]
for a choice of positive roots $\Sigma_{+}$ and 
an element $c \in \fh \wedge \fh$. It defines a 
(holomorphic) Poisson structure $\piG$ on $G$ by
\[
\piG (g) \, = \, R_g r_0 \, - \, L_g r_0,
\]
where $R_g$ and $L_g$ are respectively the right and left 
translations on $G$ by $g \in G$,
making $(G, \piG)$ into a Poisson Lie group. 
This Poisson structure is the semi-classical limit
of the quantum group corresponding to $G$ \cite{dr:bigbra} 
\cite{dr:quantum}. 
A Poisson structure on $G/H$ is said to be 
$(G, \piG)$-homogeneous
if the action map $G \times (G/H) \rightarrow G/H$ is a Poisson map
\cite{dr:homog}. 

The first result of this paper, Theorem \ref{thm_main}, 
is on the construction of a
surjective map from the set of all classical
dynamical $r$-matrices for the pair $(\fg, \fh)$ together with their domains
 to the set
of all (holomorphic) $(G, \piG)$-homogeneous 
Poisson structures on
$G/H$. More precisely, for any classical dynamical $r$-matrix
$r$ and $\lambda \in \fh^*$ such that $r(\lambda)$ is defined,
we show that the bi-vector field $\tilde{\pi}_{r(\lambda)}$
on $G$ defined by
\[
\tilde{\pi}_{r(\lambda)} \, = \, R_{g} r_0 \, - \, L_g r(\lambda)
\]
projects to a holomorphic $(G, \piG)$-homogeneous 
Poisson structure
on $G/H$ under the projection $G \rightarrow G/H$, and that all
$(G, \piG)$-homogeneous Poisson structures on $G/H$
arise this way. 
See also \cite{x-l:homog} for another interpretation
of classical dynamical $r$-matrices.

Let $K \subset G$ be a compact real form of $G$, and let $T = K \cap H$
be the maximal torus of $K$. Then $K$ also carries a natural Poisson
structure $\piK$ such that $(K, \piK)$ is a Poisson Lie group.
Theorem \ref{thm_main} is then modified to Theorem \ref{thm_compact}
which states that classical dynamical $r$-matrices give rise to 
 $(K, \piK)$-homogeneous Poisson structure on $K/T$
and that all $(K, \piK)$-homogeneous Poisson structures on $K/T$
arise this way. 

We point out that a classification of all
$(G, \piG)$ or $(K, \piK)$-homogeneous 
Poisson structures,
not necessarily on $G/H$ or on $K/T$, has already been obtained 
by E. Karolinsky
\cite{ka:homog-compact} \cite{ka:homog-complex}. We want to emphasize 
that what is brought out here is the connection of such Poisson spaces with
the CDYBE.

Among all $(K, \piK)$-homogeneous Poisson structures on $K/T$,
we single out a family denoted by $\pix$, where
$X$ is any subset of the set $S(\Sigma_{+})$
of all simple roots, $X_1 \subset X$, 
 and $\lambda \in \fh$ satisfies some
regularity condition (Theorem \ref{thm_pix}). This family exhausts all
$(K, \piK)$-homogeneous Poisson structures on $K/T$ up to
$K$-equivariant isomorphisms. Moreover, these Poisson structures
are related to each other by taking various limits of the parameter 
$\lambda$ (see Section \ref{sec_limits}).
We study several Poisson geometrical properties of this family:

The Lagrangian subalgebra of $\fg$ corresponding
to each $\pix$ is described in Section \ref{sec_lagrangian-compact}.

In Section \ref{sec_geom-X-whole}, we recall  the construction 
in \cite{e-l:Lagrangian} of
a Poisson structure $\Pi$ on the variety $\Lagr$ of all Lagrangian 
subalgebras in $\fg$ and the fact that each $(K/T, \pix)$
sits inside $(\Lagr, \Pi)$ as a Poisson submanifold (possibly up to
a covering map). The two special cases of $\pix$ when 
$X= X_1 = \emptyset$ and when $X = S(\Sigma_{+}), X_1 = \emptyset$
are considered in more detail here.

In Section \ref{sec_induction}, we show that each $\pix$ on $K/T$ can be 
obtained via Poisson induction from a Poisson structure on a smaller
manifold.

 
In Section \ref{sec_leaves-1}, we describe
the symplectic leaves
of $\pix$ when $X_1$ is the empty set.
We show that in this case $\pix$ has a finite number of symplectic leaves.
For an arbitrary $\pix$, we show that it always has at least one open
symplectic leaf.

In Section \ref{sec_modular}, we show that with
respect to a $K$-invariant volume form $\mu_0$ on $K/T$, all
the Poisson structures $\pix$ have the same modular
vector field. In the case when $X_1$ is the empty set,
we also describe the moment map for the $T$-action on
each symplectic leaf of $\qix$. 


Some applications of results in this paper are given in \cite{e-l:harm}, 
where
a Poisson
geometrical interpretation of the Kostant harmonic forms
on $K/T$ \cite{ko:63} is given using the Bruhat Poisson structure
$\pinf := \pix$ for $X = X_1 = \emptyset$.  Set $\pi_\lambda = 
\pix$ when $X = S(\Sigma_{+})$ and $X_1 = \emptyset$.
The fact that $\pi_{\lambda} \rightarrow \pi_{\infty}$ as 
$\lambda \rightarrow \infty$ is used in 
\cite{e-l:harm} to show that
the Kostant harmonic forms are limits of the usual Hodge harmonic forms.

Results in this paper also motivate our work in \cite{e-l:Lagrangian},
where, among other things, we show that there is a Poisson manifold 
$({\cal L}_0, \Pi)$ such that every $(K/T, \pix)$
is a Poisson submanifold (possibly up to a covering map)
of $({\cal L}_0, \Pi)$. In fact, ${\cal L}_0$ is an irreducible
component of the variety ${\cal L}$
of all Lagrangian subalgebras of $\fg$, and the Poisson
structure $\Pi$ is defined on all of ${\cal L}$. We show in
\cite{e-l:Lagrangian} that all the $K$-orbits in ${\cal L}$
with respect to the Adjoint action are $(K, \piK)$-homogeneous
Poisson spaces, and that every $(K, \piK)$-homogeneous
Poisson space maps to $({{\cal L}}, \Pi)$ by a 
Poisson map. Thus, $({{\cal L}}, \Pi)$ is a setting
for studying all $(K, \piK)$-homogeneous Poisson spaces.

We point out that many more properties of the Poisson structures $\pix$
can be studied, among these their Poisson cohomology, their
Poisson harmonic forms \cite{e-l:harm}, and their symplectic groupoids.
We hope to do this in the future.


{\bf Acknowledgement} The author would like to thank P. Etingof
for explaining to her the results in \cite{e-v:cdyb}
and Professors V. Drinfeld, S. Evens,  Y. Kosmann-Schwarzbach, A. Weinstein 
and P. Xu for helpful discussions. She would also like
to thank the Mathematics Department of
 Hong Kong University of Sciences and
Technology for it hospitality.

\section{The Classical Dynamical Yang-Baxter Equation}
\label{sec_cdybe}

\begin{dfn} 
\label{dfn_r} \cite{fl:cdyb} \cite{e-v:cdyb} 
{\em
A meromorphic function $ r:  \fh^* \rightarrow \fg \ot \fg$
is called a {\it classical (quasi-triangular) dynamical 
$r$-matrix for the pair}
$(\fg, \fh)$ if it satisfies the following three conditions:

1. {\it The zero weight condition:} $ ad_x r(\lambda)  =  0$
for all $x \in \fh$ and $\lambda \in \fh^*$ such that $r(\lambda)$
is defined;

2. {\it The generalized unitarity condition:}
$r^{12} + r^{21} = \varepsilon \Omega$
for some complex number $\epe$ and for all $\lambda \in \fh^*$ such that
$r(\lambda)$ is defined, 
where $\Omega \in (S^2 \fg)^{\frak g}$ is the element corresponding to
the Killing form on $\fg$;

3. {\it The Classical Dynamical Yang-Baxter Equation (CDYBE):}
\[
{\rm Alt} (d r) \, + \, [r^{12},r^{13}]\,+ \, [r^{12},r^{23}]\,
+\,[r^{13},r^{23}] \, = \, 0 \,,
\]
where, for $r = \sum_i u_i \ot v_i$, we have 
$r^{12} = \sum_i u_i \ot v_i \ot 1, \,
r^{13} = \sum_i u_i \ot 1 \ot v_i, \,
r^{23} = \sum_{i} 1 \ot u_i \ot v_i,$ 
\beqa
{\rm CYB}(r) & := & [r^{12},r^{13}] + [r^{12}, r^{23}] + [r^{13}, r^{23}] \\
& = &\sum_{i,j} [u_i, u_j] \ot v_i \ot v_j + 
u_i \ot [v_i, u_j] \ot v_j + 
u_i \ot u_j \ot [v_i, v_j],
\eeqa
and ${\rm Alt} (d r)(\lambda) \in \wedge^3 \fg$ is the 
skew-symmetrization of $dr(\lambda) \in \fh \ot \fg \ot \fg \subset
\fg \ot \fg \ot \fg$.
The complex number $\varepsilon$ is called the {\it coupling constant}
for $r$.
}
\end{dfn}

We now recall the classification of classical dynamical
$r$-matrices for the pair $(\fg, \fh)$ as given in \cite{e-v:cdyb}.
Let $\Sigma$ be the set of all roots for $\fg$ with respect to 
$\fh$. 
For each $\alpha \in \Sigma$, choose root
vectors $\ea$ and $\eb$ such that $\ll \ea, \eb \gg = 1$, where $\ll ~, ~ \gg$
is the Killing form on $\fg$. 

Let $\epe$ be a non-zero complex number, let $\mu \in \fh^*$, and let
$C = \sum_{i,j} C_{ij} dx_i \wedge dx_j$ be a closed meromorphic
$2$-form on $\fh^*$. Let $\Sigma_{+}$ be a choice of
positive roots, and let $X$ be a subset of the set $S(\Sigma_{+})$
of simple roots in $\Sigma_{+}$. For each 
$\alpha \in \Sigma$, define a (scalar-valued) meromorphic function
$\phi_{\alpha}$ on $\fh^*$
according to the rule: If $\alpha$ is a linear combination of simple
roots in $X$, then 
\[
\phi_{\alpha} (\lambda) ~ = ~ 
{\epe \over 2} \, \coth \, ({\epe \over 2} \, \ll \alpha, \lambda -  \mu\gg),
\]
where $\coth(x) = {\frac{e^x + e^{-x}}{e^x - e^{-x}}}$ is the hyperbolic
cotangent function;
Otherwise, set $\phi_{\alpha}(\lambda) = {\epe \over 2}$ if
$\alpha$ is positive and $\phi_{\alpha}(\lambda) = -{\epe \over 2}$
if $\alpha$ is negative.

\begin{thm}[Etingof-Varchenko \cite{e-v:cdyb}]
\label{thm_ev}

1. With the above choices of $\mu, C, \Sigma_{+},  X\subset  
S(\Sigma_{+}) $ and $\phi_{\alpha}$,  the
meromorphic function
$r: \fh^* \rightarrow \fg \ot \fg$ defined by
\begin{equation}
\label{eq_r-general}
r(\lambda) \, = \, {\epe \over 2}\, \Omega \,+ \,
\sum_{i,j} C_{ij}(\lambda) x_i \otimes x_j \,
+ \,
\sum _{\alpha \in \Sigma} \,
\phi_\alpha (\lambda )\,
\ea \otimes  \eb\,
\end{equation}
is a classical dynamical $r$-matrix with non-zero coupling constant
$\varepsilon$;

2. Every classical dynamical $r$-matrix 
with non-zero coupling constant has this form.
\end{thm}

\section{$r$-matrices and homogeneous Poisson structures on $G/H$}
\label{sec_poi-ongh}

\subsection{The main theorem}
\label{sec_main}

Let $r: \fh^* \rightarrow \fg \ot \fg$ be 
any classical dynamical $r$-matrix as
in Definition \ref{dfn_r}.
Let
\[
A_r(\lambda) \, = \, r(\lambda) \, - \, {\epe \over 2} \Omega
\]
be the skew-symmetric part of $r(\lambda)$. Using the fact that
$\Omega$ is symmetric and $ad$-invariant, one easily shows that the terms
$[\Omega^{ij}, A_{r}^{kl}]$ in the CDYBE for $r$ all cancel. 
Moreover, 
it is well-known that 
\[
[\Omega^{12},  \Omega^{13}]  +  
[\Omega^{12},  \Omega^{23}]  + 
[\Omega^{13},  \Omega^{23}]   =  
[\Omega^{12},  \Omega^{13}] 
 =  [\Omega^{13},  \Omega^{23}] =  -\,[\Omega^{12}, 
 \Omega^{23}]  \in  (\wedge^3 \fg)^{\frak g}.
\]
Therefore, 
$A_r$ satisfies the following modified CDYBE (see also \cite{e-v:cdyb}):
\begin{equation}
\label{eq_A1}
{\rm Alt} (d A_r) \, + \, [A_{r}^{12},A_{r}^{13}]\,+ \, 
[A_{r}^{12},A_{r}^{23}]\,
+\,[A_{r}^{13},A_{r}^{23}] \, = 
\,  {\epe^2 \over 4}\,
[\Omega^{12},\Omega^{23}]\, \in (\wedge^3 \fg)^{\frak g}.
\end{equation}
Recall that there is the Schouten bracket $[ \hspace{.2in}]$
on $\wedge \fg$. 
For $x_1, x_2, ..., x_k \in \fg$, we use the convention
\[
x_1 \wedge x_2 \wedge \cdots \wedge x_k \, = \, 
\sum_{\sigma \in S_k} {\rm sign}(\sigma) x_{\sigma(1)} \ot
x_{\sigma(2)} \ot \cdots \ot x_{\sigma(k)} \, \in \fg^{\otimes k}.
\]
Then  for $X \in \wedge^2 \fg$, the element
${\rm CYB}(X)$ and the Schouten bracket $[X, \, X]$ 
are related by \cite{dr:quantum} 
\[
{\rm CYB}(X) \, = \, [X^{12}, \, X^{13}] \, + \, [X^{12}, \, X^{23}] \, + \,
[X^{13}, \, X^{23}] \, = \, {\frac{1}{2}} \,[X, \, X].
\]
Thus, we can rewrite Equation (\ref{eq_A1}) as
\begin{equation}
\label{eq_A2}
[A_r(\lambda), \, A_r(\lambda)] \, = \, {\epe^2 \over 2}\,
[\Omega^{12},\Omega^{23}]\, - \, 2{\rm Alt}(dA_r)(\lambda).
\end{equation}
It is this form of the CDYBE that we will use to define Poisson structures on
$G/H$.

Recall \cite{dr:quantum} that a {\it classical quasi-triangular $r$-matrix with 
coupling constant $\epe$} is an element $r_0 \in \fg \ot \fg$ such that
\beqa
& & r_{0} \, + \, r_{0}^{21} \, = \, \epe \Omega\\
& & {\rm CYB} (r_0) \, = \, 0.
\eeqa

\begin{rem}
\label{rem_other-r}
{\em
If $r_0$ has the zero-weight property, i.e., if $r_0 \in
(\fg \ot \fg)^{\frak h}$, then by
Theorem \ref{thm_ev}, it must be of the form
\begin{equation}
\label{eq_constant}
r_0 \, = \, {\frac{\epe}{2}} \Omega \, + \, \sum_{i,j} c_{ij} x_i
\wedge x_j \, + \,  {\frac{\epe}{2}} \sum_{\alpha \in \Sigma_{+}}
E_{\alpha} \wedge E_{-\alpha}
\end{equation}
for some choice $\Sigma_{+}$ of positive roots and $\sum_{i,j} c_{ij} \in
\fh \wedge \fh$.
But not every quasi-triangular $r_0$ has the zero-weight property.
For example, for $\fg = \fsl(3, \C)$, we can take 
$r_0 = {\epe \over 2} (\Omega + h \wedge (e + f))$ where
$h, e$ and $f$ are the three generators with 
Lie brackets: $[h, e] = 2e, \, [h, f] = -2 f$ and
$[e, f] = h$. See \cite{b-d:r} for more examples.
}
\end{rem}

Let $r_0$ be a classical quasi-triangular $r$-matrix with coupling constant
$\epe$ (not necessarily of zero weight for $\fh$). 
Let $\Lambda = r_0 - {\frac{\epe}{2}} \Omega \in \fg \wedge \fg$
be the skew-symmetric part of $r_0$. Then, as a special case of
(\ref{eq_A2}),
$\Lambda$ satisfies
the modified Classical Yang-Baxter Equation (CYBE)
\begin{equation}
\label{eq_Lambda}
[\Lambda, \, \Lambda] \, = \, {\frac{\epe^2}{2}} [\Omega^{12}, \, \Omega^{23}].
\end{equation}
It is well known 
that the bi-vector field
$\piG$ on 
the group $G$ defined by
\begin{equation}
\label{eq_pi-on-G}
\piG (g) \, = \, R_g \Lambda \, - \, L_g \Lambda,
\end{equation}
where for $R_g$ and $L_g$ denote respectively the right and left translations
from the identity element to $g$, defines a Poisson structure on $G$, and that
$(G, \piG)$ is a Poisson Lie group \cite{dr:quantum} \cite{sts:rmatr}.

\begin{rem}
\label{rem_holom}
{\em
The meaning of the terms $R_g \Lambda$ and
$L_g \Lambda$ needs further explanation. Denote by $J$ the 
complex structure on $\fg$ induced by that on $G$. Then we can
identify $(\fg, J)$ with the holomorphic tangent space
$T_{e}^{1,0}G$ of $G$ at $e$ via
$\fg \ni x \mapsto  {\frac{1}{2}} (x  -  i J(x)).$
For $\Lambda \in \fg \wedge \fg$, we regard $\Lambda$ as an
element in
$\wedge^2 T_{e}^{1,0}G$. Then, $L_g \Lambda$ (resp. 
$R_g \Lambda$), for $g \in G$,
is understood to be the image in 
$\wedge^2 T_{g}^{1,0}G$ of $\Lambda$ by the left (resp. right)
 translation by $g$. Thus the bi-vector field $\piG$
on $G$ in (\ref{eq_pi-on-G}) is holomorphic.
All Poisson structures in this section are assumed to be holomorphic.
}
\end{rem}

Recall that an action of 
the Poisson Lie group $(G, \piG)$ 
on a Poisson manifold $P$ is said to
be Poisson if the action map 
$G \times P \rightarrow P:  (g, \, p) \mapsto gp $
is a Poisson map, where $G \times P$ is equipped with the
product Poisson structure. When the action of $G$ on $P$ is transitive, 
the Poisson structure on $P$ is said to be $(G, \piG)$-homogeneous
\cite{dr:homog}.
The following theorem makes a connection between classical dynamical
$r$-matrices and $(G, \piG)$-homogeneous Poisson structures on 
$G/H$.

\begin{thm}
\label{thm_main}
Let $r_0 = {\frac{\epe}{2}} \Omega + \Lambda$ be any 
classical quasi-triangular 
$r$-matrix (not necessarily of zero-weight) with skew-symmetric part $\Lambda$.
Let $r(\lambda) = {\frac{\epe}{2}} \Omega + A_r(\lambda)$
be any classical dynamical $r$-matrix for the pair $(\fg, \fh)$ as
in Definition \ref{dfn_r}. 
For each value $\lambda$ such that $r(\lambda)$ is defined, 
define a bi-vector field $\tlp_{r(\lambda)}$ on $G$ by
\[
\tlp_{r(\lambda)} (g) \, = \, R_g \Lambda \, - \, 
L_g A_r(\lambda), \hspace{.2in} g \in G.
\]
Let $\pi_{r(\lambda)} = p_* \tlp_{r(\lambda)}$ 
be the projection of $\tlp_{r(\lambda)}$ to $G/H$
by the map $p: G \rightarrow G/H: g \mapsto gH$. Then

1) $\pi_{r(\lambda)}$ is well-defined and it defines a Poisson
structure on $G/H$;

2) Equip $G$ with the Poisson structure $\piG$ as defined by
(\ref{eq_pi-on-G}). Then 
$\pi_{r(\lambda)}$ is a $(G, \, \piG)$-homogeneous 
Poisson structure on $G/H$.

3) When $r_0$ has the zero-weight property, i.e., 
$r_0 \in (\fg \ot  \fg)^{\frak h}$, every 
$(G, \piG)$-homogeneous Poisson structure on $G/H$ arises this way.
\end{thm}

The rest of this section is devoted to the proof of this theorem.
We first prove the first two parts.

\bigskip
\noindent
{\bf Proof of 1) and 2) in Theorem \ref{thm_main}.} 
It follows from $A_r(\lambda) \in (\wedge^2 \fg)^{\frak h}$
that $\pi_{r(\lambda)}$ is well-defined. To show that 
$\pi_{r(\lambda)}$ defines a Poisson structure on $G/H$, we calculate the
Schouten bracket $[\pi_{r(\lambda)}, \, \pi_{r(\lambda)}]$ of
$\pi_{r(\lambda)}$ with itself.
Set $\Lambda^R(g) =  R_g \Lambda$
and $A_r(\lambda)^L(g) = L_g A_r(\lambda)$. Then 
$\tlp_{r(\lambda)} = \Lambda^R - A_r(\lambda)^L$. Hence
\beqa
[\tlp_{r(\lambda)}, \, \tlp_{r(\lambda)}] & = & [\Lambda^R, \, \Lambda^R] 
\, - \, 2 [\Lambda^R, \, A_r(\lambda)^L] \, + \, [A_r(\lambda)^L, \, 
A_r(\lambda)^L]\\
& = & -[\Lambda, \, \Lambda]^R \, + \, [A_r(\lambda), \, A_r(\lambda)]^L\\
& = & -2{\rm Alt}(dA_r(\lambda))^{L} \in (\fh \wedge \fg \wedge \fg)^{L},
\eeqa
where in the last step, we used Equations (\ref{eq_A2}) 
and (\ref{eq_Lambda}).
This shows that $\tlp_{r(\lambda)}$ is in general not a Poisson bi-vector
 field on $G$. 
However, for $\pi_{r(\lambda)} = p_* \tlp_{r(\lambda)}$, we have
\[
[\pi_{r(\lambda)}, \, \pi_{r(\lambda)}]  =  p_* [\tlp_{r(\lambda)}, \, 
\tlp_{r(\lambda)}]\\
 =  - 2p_* {\rm Alt}(dA_r(\lambda))^L = 0.
\]
Therefore, $\pi_{r(\lambda)}$ is a  Poisson structure on $G/H$.
Now for any $g_1$ and $g_2 \in G$, we have
\beqa
\tlp_{r(\lambda)}(g_1 g_2) & = & R_{g_1 g_2} \Lambda \, - \, 
L_{g_1 g_2} A_r(\lambda) \\
& = & L_{g_1} (R_{g_2} \Lambda  \, - \, L_{g_2} A_r(\lambda))
\, + \, R_{g_2} (R_{g_1} \Lambda - L_{g_1} \Lambda)\\
& = & L_{g_1} \tlp_{r(\lambda)} ( g_2) \, + \, R_{g_2} \piG(g_1).
\eeqa
Projecting $\tlp_{r(\lambda)}$ to $\pi_{r(\lambda)}$, 
this says that the action map of
$G$ on $G/H$ by left translations is a Poisson map. Thus  
$\pi_{r(\lambda)}$ is a $(G, \, \piG)$-homogeneous Poisson 
structure
on $G/H$. 
This finishes the proof of 1) and 2) in Theorem \ref{thm_main}.

\bigskip
We now prove 3) of Theorem \ref{thm_main}. 
Assume that $r_{0} \in (\fg \ot \fg)^{\frak h}$. Then
by Theorem \ref{thm_ev}, it must be of the form (\ref{eq_constant})
for some choice $\Sigma_{+}$ of positive roots and some
$\sum_{i,j} u_{ij} x_i \wedge x_j \in \fh \wedge \fh$.
 
Let $e = eH $ be the base point of $G/H$. 
Recall \cite{dr:homog} that
a $(G, \piG)$-homogeneous Poisson structure
$\pi$ on $G/H$ is determined by its value $\pi(e)$ at $e$
in such a way that 
\begin{equation}
\label{eq_eo}
\pi(gH) ~ = ~ L_g \pi(e) ~ + ~ p_* \piG(g).
\end{equation} 
Moreover, since $\piG(g) = 0$ for $g \in H$ (this is why we need 
the zero weight condition on $r_0$), we see that $\pi(e)$ is
$H$-invariant, i.e.,
\[
\pi (e) ~ \in  \wedge^2 T_{e} (G/H)^{H} ~ \cong ~
(\wedge^2 (\fg / \fh))^{H}. 
\]

Let $\fn_{+}$ and $\fn_{-}$ be the nilpotent Lie subalgebras
of $\fg$ spanned by the root vectors for the
 roots in $\Sigma_{+}$ and $-\Sigma_{+}$
respectively.
Identify $\fg/ \fh \cong \fn_{-} + \fn_{+}$.
 
\begin{lem}
\label{lem_three}
Write
\begin{equation}
\label{eq_pi-eo}
\pi(e) = \sum_{\alpha  \in \Sigma_{+}} ({\frac{\epe}{2}} - \phi_{\alpha})
\ea \wedge \eb ~ \in (\wedge^2 (\fg /\fh))^{H}
\end{equation}
and set $\phi_{-\alpha} = - \phi_{\alpha}$.
Then the bi-vector field $\pi$ on $G/H$ defined by (\ref{eq_eo})
is Poisson if and only if the function
$\phi: \Sigma \rightarrow \C$ satisfies
\begin{equation}
\label{eq_phi}
\phi_{\alpha} \phi_{\beta} + \phi_{\beta} \phi_{\gamma}
+ \phi_{\gamma} \phi_{\alpha}\,  = \, 
- {\frac{\epe^2}{4}}, ~~  {\em whenever} ~~
\alpha, \beta, \gamma \in \Sigma ~ {\rm and} ~ \alpha + \beta + \gamma =0.
\end{equation}
\end{lem}
 
\noindent
{\bf Proof of Lemma \ref{lem_three}.} 
For any given $\pi(e)$ in the form of 
(\ref{eq_pi-eo}), set
\[
A \, = \, \sum_{\alpha \in \Sigma_{+}} 
\phi_{\alpha} \ea \wedge \eb \in \wedge^2 \fg
\]
and introduce the following bi-vector field $\hlp$ on $G$:
\[
\hlp (g) \, = \, R_g \Lambda \, - \, L_g A.
\]
Then $\pi = p_{*} \hlp$, and hence $[\pi, \pi] = p_{*}[\hlp, \hlp]$.
But as in the proof of 1) of Theorem \ref{thm_main}, we have
\[
[\hlp, \hlp] \, = \, [\Lambda^R, \, \Lambda^R]
\, - \, 2 [\Lambda^R, \, A^L] \, + \, [A^L, \, A^L] \, 
 =\,  - \, [\Lambda, \, \Lambda]^R \, + \, [A, \, A]^L.
\]
Since $\Lambda$ satisfies the modified CYBE (\ref{eq_Lambda}), by writing
\[
B \, = \, [A, \, A] \, - \, {\frac{\epe^2}{2}} 
[\Omega^{12}, \, \Omega^{23}] \, \in \wedge^3 \fg,
\]
we see that $[\hlp, \hlp] = B^L$, the left invariant $3$-vector field
on $G$ with value $B$ at $e$.  Thus $[\pi, \pi] = 0$ if and only if
$B \in \fh \wedge \fg \wedge \fg$, or, if and only if
\[
[A, \, A] \, = \, {\frac{\epe^2}{2}}  
[\Omega^{12}, \, \Omega^{23}]  \, \, \,{\rm mod} \, \,\fh \wedge \fg 
\wedge \fg.
\]
A direct calculation shows that
\beqa
[A, \, A] &  = & \sum_{\alpha \in \Sigma} \phi_{\alpha}^{2} h_{\alpha}
\wedge \ea \wedge \eb \\
& & \, \, -2  \sum_{[(\alpha, \beta, \gamma)] \in \tilde{\Sigma}^3} 
(\phi_{\alpha} \phi_{\beta} + \phi_{\beta} \phi_{\gamma} + 
\phi_{\gamma} \phi_{\alpha}) N_{\alpha, \beta} \ea \wedge E_{\beta} \wedge
E_{\gamma}
\eeqa
and
\[
[\Omega^{12}, \, \Omega^{23}] \,= \, {\frac{1}{2}} \sum_{\alpha \in \Sigma}
 h_{\alpha} 
\wedge \ea \wedge \eb 
 \, +\, \sum_{[(\alpha, \beta, \gamma)] \in \tilde{\Sigma}^3} 
 N_{\alpha, \beta} \ea \wedge E_{\beta} \wedge 
E_{\gamma},
\]
where $h_{\alpha} = [\ea, \eb] \in \fh, \, [E_{\alpha}, E_{\beta}] = 
N_{\alpha, \beta} E_{\alpha + \beta}$ when $\alpha, \beta \in \Sigma$ 
and $\alpha + \beta \in \Sigma$, and 
the summation over
$[(\alpha, \beta, \gamma)] \in \tilde{\Sigma}^3$ means that 
the summation index runs over all triples
$(\alpha, \beta, \gamma) \in \Sigma^3$ such that $\alpha + 
\beta + \gamma = 0$ but two such triples are considered the same
if they only differ by a reordering of the three roots.
It then follows immediately that $\pi$ is a Poisson structure on 
$G/H$ if and 
only if Condition (\ref{eq_phi}) is satisfied. This finishes the  proof of
Lemma \ref{lem_three}.

\bigskip
It now remains to classify all odd functions $\phi$ on $\Sigma$
such that Condition (\ref{eq_phi}) is satisfied.  Note that the Weyl
group $W$ for $(\fg, \fh)$ acts on the set
of such functions by $(w \cdot \phi)_{\alpha} := \phi_{w \alpha}$.
We say that two such functions $\phi$ and $\psi$ are $W$-related
if $\psi = w\cdot \phi$ for some $w \in W$.

\begin{nota}
\label{nota_X}
{\em
Let $S(\Sigma_{+})$ be the set of simple roots in $\Sigma_{+}$. 
For a subset $X$ of $S(\Sigma_{+})$, we will use $[X]$ to denote
the set of roots in $\Sigma$ that
are in the linear span of $X$.
Also set
\[
\fh_{\tx} \, = \, {\rm span}_{{\Bbb C}} \{h_{\gamma} = [E_{\gamma}, 
E_{-\gamma}]: \,
\gamma \in X\}.
\]
}
\end{nota}

\begin{lem}
\label{lem_alcove}
For any $X \subset S(\Sigma_{+})$ and $h \in \fh_{\tx}$ such that
$\alpha(h) \notin \pi i \Z$ for any $\alpha \in
[X]$, where $\pi = 3.14159...$ (we hope that there is no confusion
between this notation of $\pi = 3.14159...$ and $\pi$ as 
a Poisson structure),
 and $\Z$ is the set of integers, 
define $\phi: \Sigma \rightarrow
\C$ by
\[
\phi_{\alpha} \, = \, \left\{ \begin{array}{ll}
{\frac{\epe}{2}} \coth  \alpha(h), 
 & \, \, \alpha \in  
[X]\\ \vspace{-.05in}& \vspace{-.05in} \\
 {\epe \over 2}, &  \, \, \alpha \in \Sigma_{+} \backslash
[X]\\ \vspace{-.05in}& \vspace{-.05in} \\
-{\epe \over 2}, & \, \,  \alpha \in -( \Sigma_{+} \backslash 
[X]).
\end{array} \right.
\]
Then

(1) $\phi$ satisfies Condition (\ref{eq_phi});

(2) Any odd function $\phi: \Sigma \rightarrow \C$ satisfying Condition
(\ref{eq_phi}) is $W$-related to one obtained this way.
\end{lem}

\noindent
{\bf Proof.} (1) can be checked directly. We only show (2).
Suppose that  $\phi: \Sigma \rightarrow \C$ satisfies Condition
(\ref{eq_phi}). Set $Y = \{\alpha \in \Sigma: \phi_{\alpha} = 
{\epe \over 2} \}$. Then because of (\ref{eq_phi}),
$Y$ has two properties:

(A). If $\alpha, \beta \in Y$ and $\alpha + \beta \in \Sigma$, then $\alpha +
\beta \in Y$;

(B). If $\alpha \in Y$, then $-\alpha \not\in Y$.
 
\noindent
It follows \cite{e-v:cdyb} that there
exists a choice of positive roots
$\Sigma_{+}^{'}$ such that $Y \subset \Sigma_{+}^{'}$.
Since there exists  $w \in W$ such that 
$w \Sigma_{+}^{'} = \Sigma_{+}$, by considering
$w \cdot \phi$ instead of $\phi$, we can assume that
$\Sigma_{+}^{'} = \Sigma_{+}$. Set $X = S(\Sigma_{+})
\cap (\Sigma_{+} \backslash Y)$. Since  Condition (\ref{eq_phi})
implies that $Y$ has the additional property:

(C) If $\alpha \in Y, \beta \in \Sigma \backslash (-Y)$ are such that
$\alpha + \beta \in \Sigma$, then $\alpha + \beta \in Y$, 

\noindent
we claim that
$\Sigma_{+} = ([X] \cap \Sigma_{+})  \cup Y$ is a disjoint union. Indeed,
suppose that $\alpha \in [X] \cap \Sigma_+$. We first use induction
on the height ${\rm ht}(\alpha)$
 of $\alpha$ with respect to $S(\Sigma_{+})$
 to show that $\alpha \notin Y$.
If ${\rm ht}(\alpha) = 1$, then $\alpha$ is simple, so $\alpha \notin Y$ 
by definition. Suppose that ${\rm ht}(\alpha) = k$.
We can \cite{sr:lie} write
$\alpha$ as $\alpha = \alpha_1 + \cdots + \alpha_k$ such that
each $\alpha_j$ is in $X$ and that each $\alpha_1 +
 \cdots + \alpha_j$ is a root, for $j = 1, ..., k$.
Set $\alpha^{'} = \alpha_1 + \cdots + \alpha_{k-1}$. By induction
assumption, $\alpha^{'} \notin Y$. If $\alpha \in Y$, then
we know by (C) that $\alpha_k = \alpha - \alpha^{'} \in Y$ which is a
contradiction. Thus $\alpha \notin Y$. This shows that
$([X] \cap \Sigma_+) \cap Y = \emptyset$.
 Next, suppose that $\alpha \in \Sigma_{+} \backslash Y$. We use induction on
${\rm ht}(\alpha)$ again to show that $\alpha \in [X]$.
If ${\rm ht}(\alpha) = 1$, then $\alpha \in X \subset
[X]$ by the definition of $X$. Suppose that ${\rm ht}(\alpha) = k$.
Write $\alpha$ as
$\alpha = \alpha^{'} +  \alpha_k$, where $\alpha^{'} \in \Sigma_{+}$
and $\alpha_k$ is a simple
root.
If $\alpha_k \in Y$. Then by (C), we have $-\alpha^{'} = \alpha_k -\alpha
\in Y$ which is absurd. Thus $\alpha_k \notin Y$, so $\alpha_k \in X$.
If $\alpha^{'} \in Y$, then
again by (C), we have $-\alpha_k = \alpha^{'} - \alpha \in Y$ which is
also absurd, so $\alpha^{'} \notin Y$. By induction assumption, 
$\alpha^{'} \in [X]$. Thus $\alpha \in [X]$.
Hence we have shown that $\Sigma_{+} = ([X] \cap \Sigma_+) \cup Y$ is
a disjoint union.

For $\gamma \in X$, since $\phi_{\gamma} \neq \pm {\epe \over 2}$,
there exists $\lambda_{\gamma} \in \C, \lambda_{\gamma} \notin
\pi i \Z,$  such that $\phi_{\gamma} = 
{\epe \over 2} \coth \lambda_{\gamma}$.
Choose $h \in \fh_{\tx}$ such that $\gamma(h) = \lambda_\gamma$
for every $\gamma \in X$.
We now show that
$\alpha(h) \notin \pi i \Z$ and that
$\phi_{\alpha} = {\epe \over 2} \coth \alpha(h)$
 for all $\alpha \in [X] \cap \Sigma_+$
by using induction on the height ${\rm ht}(\alpha)$.
This is true when ${\rm ht} (\alpha) = 1$. Suppose that
${\rm ht}(\alpha) = k$. As before, write 
$\alpha = \alpha^{'} + \alpha_k$, where $\alpha^{'} \in [X] \cap \Sigma_+, 
{\rm ht}(\alpha^{'}) = k-1$, and $\alpha_k \in X$. Then by induction assumption,
$\alpha^{'}(h) \notin \pi i \Z$ and $\phi_{\alpha^{'}} = 
{\epe \over 2} \coth \alpha^{'}(h)$. By Condition 
(\ref{eq_phi}),
\[
-\phi_{\alpha} (\phi_{\alpha^{'}} + \phi_{\alpha_k}) \, = \,
-{\epe^2 \over 4} - \phi_{\alpha^{'}} \phi_{\alpha_k}.
\]
If $\phi_{\alpha^{'}} + \phi_{\alpha_k} = 0$, we would have 
$\phi_{\alpha^{'}} \phi_{\alpha_k} = -{\epe^2 \over 4}$ and thus
$\phi_{\alpha^{'}} = \pm {\epe \over 2}$ and
$\phi_{\alpha_k} = \mp {\epe \over 2}$. This is not possible
since $([X] \cap \Sigma_+)\cap Y = \emptyset$. Thus $\phi_{\alpha^{'}} +
\phi_{\alpha_k} \neq 0$, so $\alpha(h) = \alpha^{'}(h) + \alpha_k(h)
 \notin \pi i \Z$, and
\[
\phi_{\alpha} \, = \, {{\epe^2 \over 4} + \phi_{\alpha^{'}} \phi_{\alpha_k}
\over \phi_{\alpha^{'}} + \phi_{\alpha_k}} \, = \, 
 {\epe \over 2} \coth \alpha(h).
\]
\qed

We now continue with the proof of (3) of Theorem 
\ref{thm_main}. Let $\pi$ be a $(G, \piG)$-homogeneous
Poisson structure on $G/H$. Then by Lemmas \ref{lem_three} and
\ref{lem_alcove}, there exist a choice $\Sigma_{+}^{'}$
of positive roots, a subset $X^{'}$ of
the set of simple roots in $\Sigma_{+}^{'}$, and an element 
$\lambda_0 \in \fh^*$ such that 
$\pi = \pi_{r_{X^{'}}(\lambda_0)}$, where
\begin{equation}
\label{eq_x-prime}
r_{\tx^{'}}(\lambda) \, = \, {\epe \over 2} \Omega\, +\, {\epe \over 2}
\sum_{\alpha \in [X^{'}] \cap \Sigma_{+}} 
\coth {\epe \over 2} \ll \alpha, \lambda \gg
\ea \wedge \eb \, + \, {\epe \over 2} 
\sum_{\alpha \in \Sigma_{+}^{'} \backslash [X^{'}]}
\ea \wedge \eb
\end{equation}
is a classical dynamical $r$-matrix for the pair
$(\fg, \fh)$.
This proves part (3) of Theorem \ref{thm_main}.
\qed

\subsection{The Poisson structures $\pi_{r_{\tx}(\lambda)}$ on $G/H$}
\label{sec_onGH}

In this section, we consider in more detail
the case when the Poisson structure on $G$ is defined by
a classical
quasi-triangular $r$-matrices $r_0$
{\it with the zero weight property}. In other words,
we fix a choice $\Sigma_{+}$ of positive roots,
and consider $r_0$ of the form
\begin{equation}
\label{eq_r0-special}
r_0 \, = \, {\frac{\epe}{2}} \Omega \, + \,
\sum_{i,j} c_{ij} x_i \wedge x_j \, + \,
{\frac{\epe}{2}} \sum_{\alpha \in \Sigma_{+}} \ea \wedge \eb,
\end{equation}
where $\sum_{i,j} c_{ij} x_i \wedge x_j \in \fh \wedge \fh$.
When $\sum_{i,j} c_{ij} x_i \wedge x_j = 0$, the corresponding
$r_0$ is often called the
standard $r$-matrix. The corresponding
Poisson structure $\piG$ on $G$ is the semi-classical
limit of the quantum group corresponding to $G$ \cite{dr:quantum}.

For $X \subset S(\Sigma_{+})$, set
\begin{equation}
\label{eq_rx}
r_{\tx} (\lambda) \, = \, {\frac{\epe}{2}} \Omega \, + \,
{\frac{\epe}{2}} \sum_{\alpha \in [X] \cap \Sigma_{+}} 
\coth {\frac{\epe}{2}}
\ll \alpha, \lambda \gg \ea \wedge \eb \, + \, {\frac{\epe}{2}}
\sum_{\alpha \in \Sigma_{+} \backslash [X] } \ea \wedge \eb.
\end{equation}
Clearly, the domain $D(r_{\tx})$ of $r_{\tx}$
consists of those $\lambda \in \fh^*$
such that $\ll \lambda, \alpha \gg \notin {2 \pi i \Z \over \epe}$
for all $\alpha \in [X]$. For each such $\lambda$, we have
the $(G, \piG)$-homogeneous Poisson
structure $\pi_{r_{\tx}(\lambda)}$ on $G/H$: let $p_* \piG$
be the projection to $G/H$ of $\piG$ by
$p: G \rightarrow G/H: g \mapsto gH$. Then
\[
\pi_{r_{\tx}(\lambda)} \, = \, p_* \piG \, + \, 
\left( \sum_{\alpha \in [X] \cap \Sigma_{+}} 
{\epe \over 1 - e^{\varepsilon \ll \a, \lambda \gg}}
\ea \wedge \eb \right)^L,
\]
where where the second 
term  on the right hand side
is the $G$-invariant bi-vector field on $G/H$
whose value at $e = eH$ is the expression given in the parenthesis.

\begin{thm}
\label{thm_G}
With the Poisson structure $\piG$ on $G$ defined by $r_0$ in
(\ref{eq_r0-special}), every
holomorphic $(G, \piG)$-homogeneous Poisson structure
on $G/H$ is isomorphic, via a $G$-equivariant diffeomorphism,
to a $\pi_{r_{\tx}(\lambda)}$ for some
subset $X \subset S(\Sigma_{+})$ and $\lambda \in D(r_{\tx})$,
where $r_{\tx}$ is given in (\ref{eq_rx}).
\end{thm}
 
\noindent
{\bf Proof.} Let $\pi$ be a $(G, \piG)$-homogeneous Poisson structure
on $G/H$. By Theorem \ref{thm_main}, we know that there
exists a choice $\Sigma_{+}^{'}$ of positive roots
and a subset $X^{'}$ of the set of simple roots in $\Sigma_{+}^{'}$ such
that $\pi = \pi_{r_{\tx^{'}}(\lambda_0)}$ for some $\lambda_0 \in \fh^*$,
where $r_{\tx^{'}}$ is the classical dynamical $r$-matrix given by
(\ref{eq_x-prime}). Let $\Lambda = r_0 -{\epe \over 2} \Omega$ and
 let $A_{\tx^{'}}(\lambda_0)$ be the skew-symmetric part of 
$r_{\tx^{'}}(\lambda_0)$. Then recall from Section \ref{sec_poi-ongh}
that $\pi = p_* \hat{\pi}^{'}$,
where 
$\hat{\pi}^{'}$ is the bi-vector field on $G$ given by
\[
\hat{\pi}^{'} (g) \, = \, R_g \Lambda \, - \, L_g A_{\tx^{'}}(\lambda_0),
\hspace{.2in} g \in G.
\]
Pick  $w \in W$ such that $w \Sigma_{+}^{'} = \Sigma_{+}$. Set
$X = w X^{'}$.
Let $\dot{w}$ be a representative of $w$ in $G$. We will use 
$R_{\dot{w}^{-1}}$ to denote the right translation on $G$ by $\dot{w}^{-1}$
as well as the induced diffeomorphism on $G/H$. Then for any $g \in G$,
\[
R_{\dot{w}^{-1}} \hat{\pi}^{'} (g) \, = \, 
R_{\dot{w}^{-1}g} \Lambda \, - \, L_g L_{\dot{w}^{-1}} {\rm Ad}_{\dot{w}}
A_{\tx^{'}}(\lambda_0) \, = \, R_{g \dot{w}^{-1}}
\Lambda \, - \,  L_{g \dot{w}^{-1}} A_{\tx}(w \lambda_0),
\]
where $A_{\tx}$ is the skew-symmetric part of the $r$-matrix
$r_{\tx}$ given by (\ref{eq_rx}).
It follows from the definition of $\pi_{r_{\tx}(w \lambda_0)}$
that $\pi = R_{\dot{w}} \pi_{r_{\tx}(w \lambda_0)}$. The map $R_{\dot{w}}:
G/H \rightarrow G/H$ is $G$-equivariant.
\qed

\subsection{Comparison with Karolinsky's classification}
\label{sec_karolin}

 When $\sum_{ij} c_{ij} x_i \wedge x_j = 0$ in the definition of $r_0$,
all $(G, \piG)$-homogeneous
Poisson structures on $G/H$ have been classified by Karolinsky
\cite{ka:homog-complex} by using Drinfeld's theorem on
Poisson homogeneous spaces. We now look at the Poisson structures
$\pi_{r_{\tx}(\lambda)}$ on $G/H$ in terms of Karolinsky's classification.

Recall that the double Lie algebra associated to the Poisson Lie
group $(G, \piG)$ can be identified with the direct sum Lie algebra
$\fd = \fg + \fg$ equipped with the ad-invariant non-degenerate
scalar product given by
\[
\la (x_1, x_2), \, (y_1, y_2) \ra \, = \, {1 \over \epe}
(\ll x_2, y_2 \gg \, - \, \ll x_1, y_1 \gg).
\]
The Lie algebra $\fg$ is identified with the diagonal of $\fd$, and
the
Lie algebra $\fg^*$ is identified with the subspace
\[
\fg^* \cong \{(x_{-}, x_{+}): ~ x_{\pm} \in \fb_{\pm},
\, \, (x_{-})_{\frak h} + (x_{+})_{\frak h} = 0 \}.
\]
Here, $\fb_{\pm} = \fh + \fn_{\pm}$
and $(x_{\pm})_{\frak h} \in \fh$ is the $\fh$-component of $x_{\pm}$.
A theorem of Drinfeld \cite{dr:homog} says that
$(G, \piG)$-homogeneous
Poisson structures on $G/H$ correspond to Lagrangian (with respect to
the scalar product $\la \, , \, \ra$) subalgebras $\fl$ of the double
$\fd \cong \fg + \fg$ such that $\fl \cap \fg = \fh$.
 
\begin{thm}[Karolinsky] \cite{ka:homog-complex}
\label{thm_karo}
Lagrangian subalgebras $\fl$ of $\fg + \fg$ such that
$\fl \cap \fg = \fh$ are in $1-1$ correspondence with triples
$(\fp, \fp^{'}, \eta)$, where $\fp$ and $\fp^{'}$ are parabolic
subalgebras of $\fg$ such that $\fq = \fp \cap \fp^{'}$ is
the Levi subalgebra, $\fh \subset \fq$, and
$\eta$ is an interior orthogonal automorphism of $\fq$ with
$\fq^{\eta} = \fh$. If $(\fp, \fp^{'}, \eta)$ is
such a triple, the corresponding subalgebra $\fl$ of $\fg + \fg$ is
$\fl = \{(x^{'}, x) \in \fp^{'} \times \fp: \,
\eta(x_{\frak q}^{'}) = x_{\frak q} \},$ where $x_{\frak q} \in
\fq$ (resp. $x_{\frak q}^{'} \in \fq^{'}$) is the projection of
$x$ (resp. $x^{'}$) to $\fq$ with respect to the Levi decomposition
of $\fp$ (resp. $\fp^{'}$).
\end{thm}
 
For a $(G, \piG)$-homogeneous Poisson structure 
$\pi$ on $G/H$,
the Lagrangian subalgebra $\fl_{\pi(e)}$ of $\fg + \fg$
is by definition \cite{dr:homog}
\[
\fl_{\pi(e)} \, = \, \{x +  \xi: \, x \in \fg, \, \xi \in \fg^*, \,
\xi|_{\frak h} = 0, \, {\rm and} \, \xi \backl \pi(e) = x + \fh \}.
\]
For $\pi(e)$ of the form
$\pi(e)  =  \sum_{\alpha \in \Sigma_{+}} ({\frac{\epe}{2}} -
\phi_{\alpha}) \ea \wedge \eb$,
it is an easy calculation to see that
\[
\fl_{\pi(e)} \, = \, \fh \, + \, {\rm span}_{\Bbb C} \{\xi_{\alpha}:
\, \alpha \in \Sigma \},
\]
where for $\alpha \in \Sigma$,
\[
\xi_{\alpha} \, = \, \left(\, 
(\phi_{\alpha} - {\frac{\epe}{2}}) \ea, \, \,
(\phi_{\alpha} + {\frac{\epe}{2}}) \ea \right) \, \in \, \fg + \fg.
\]
Thus, for the Poisson structure $\pi_{r_{\tx}(\lambda)}$ on $G/H$, we have
\[
\xi_{\alpha} \,  = \, \left\{ \begin{array}{ll} (- \epe \ea, \, 0)
& {\rm if} ~
\alpha \in -Y \vspace{.1in} \\
{\epe \over e^{\varepsilon \ll \alpha, \lambda \gg}-1}
(\ea, \, e^{ \epe \ll \alpha, \lambda \gg }
\ea) & {\rm if} \, \alpha \in  [X] \vspace{.1in} \\
(0, \, \epe \ea) & {\rm if} \, \alpha \in Y.
\end{array} \right.,
\]
where $Y = \Sigma_+ \backslash [X]$.
Let
\[
\fp_{\tx} \,  = \, 
\fh \, + \, {\rm span}_{\Bbb C} \{\ea: \, \alpha
\in [X] \cup Y \} 
\]
be the parabolic subalgebra of $\fg$ defined by $X$, and let
\[
\fp^{'}_{\tx} \,  = \, 
\fh \, + \, {\rm span}_{\Bbb C}
\{\ea: \, \alpha \in [X] \cup (-Y) \}
\]
be its opposite parabolic subalgebra.
Set
\begin{equation}
\label{eq_gx}
\fm_{\tx} \, = \,
 \fh \, + \, {\rm span}_{\Bbb C} \{\ea: \alpha \in [X] \}
\end{equation}
so that $\fm_{\tx} = \fp_{\tx} \cap \fp^{'}_{\tx}$. Let 
$\eta$ be the 
 interior automorphism of $\fm_{\tx}$ given by 
${\rm Ad}_{e^{\varepsilon h_{\lambda}}}$, where $h_{\lambda} \in \fh$
corresponds to $\lambda \in \fh^*$ under the Killing form.
Then the triple $(\fp_{\tx}^{'}, \fp_{\tx}, \eta)$
is the one corresponding to the Poisson structure 
$\pi_{r_{\tx}(\lambda)}$
in the Karolinsky classification.

\section{$r$-matrices and homogeneous Poisson structures on $K/T$}
\label{sec_compact}

We pick a compact real form $\fk$ of $\fg$ as follows: For each
$\alpha \in \Sigma_+$, set
\[
X_{\alpha} = \ea - \eb, \hspace{.4in}
Y_{\alpha} = i(\ea + \eb)
\]
and $h_{\alpha} = [\ea, \eb]$. Then the real subspace
\[
\fk \, = \, {\rm span}_{\Bbb R} \{ ih_{\alpha}, X_{\alpha}, Y_{\alpha}:
\alpha \in \Sigma_+ \}
\]
is a compact real form of $\fg$. Set
$\ft = {\rm span}_{\Bbb R} \{ih_{\alpha}: \alpha \in \Sigma \} \subset \fk$.
Let $K$ and $T \subset K$ be respectively the connected compact
subgroups of $G$ with Lie algebras $\fk$ and $\ft$.
 
It is well-known \cite{soi:compact} that every Poisson
structure $\piK$ on $K$ such that $(K, \piK)$ is a Poisson
Lie group is of the form
\begin{equation}
\label{eq_on-K}
\piK(k) \, = \, R_k \Lambda \, - \, L_k \Lambda,
\end{equation}
where
\begin{equation}
\label{eq_lambda-u}
\Lambda \, = \, u \, - \, { i \epe  \over 2} \sum_{ \alpha \in \Sigma_{+}}
{X_{\alpha} \wedge Y_{\alpha} \over 2} \, \in \, \fk \wedge \fk
\end{equation}
for some $u \in \ft \wedge \ft$, an imaginary complex number $\epe$ and
a choice $\Sigma_{+}$ of positive roots.
 
In this section, we will show how $(K, \piK)$-homogeneous
Poisson structures on $K/T$ are related to classical
dynamical $r$-matrices. We remark again that one classification
of all $(K, \piK)$-homogeneous Poisson spaces (by the corresponding
Lagrangian Lie subalgebras) has
been given by Karolinsky \cite{ka:homog-compact}.
 
If we regard $\wedge \fg$ as a real vector space, then
\[
\wedge \fk \lrw \wedge \fg: \, \, \wedge^l \fk \ni
x_1 \wedge \cdots \wedge x_l \Map x_1 \wedge \cdots \wedge x_l \in
\wedge^l \fg
\]
is an embedding of $\wedge \fk$ into $\wedge \fg$ as a real subspace.
This embedding also preserves the Schouten bracket. Thus,
for $A \in \wedge^2 \fk$ of the form
\[
A \, = \,\sum_{\alpha \in \Sigma_{+}} a_{\alpha} {\Xa \wedge \Ya \over 2},
\hspace{.3in} a_{\alpha} \in {\Bbb R} \, \, \, {\rm for } \, \,  \alpha 
\in \Sigma_{+},
\]
we can calculate $[A, A] \in \wedge^3 \fk$
by first writing $A = \sum_{\alpha \in \Sigma_{+}} 
ia_{\alpha} \ea \wedge \eb \in
\wedge^2 \fg$ and calculate $[A, A]$ inside $\wedge \fg$. Indeed, as in
the proof of Lemma \ref{lem_three}, in $\wedge^3 \fg$ we have
\begin{eqnarray}
\nonumber
[A, A] & = & {1 \over 2}\sum_{\alpha \in \Sigma_{+}} a_{\alpha}^{2}
(ih_{\alpha} \wedge \Xa \wedge \Ya) \\
\label{eq_in-k}
& & + \, 2\sum_{[(\alpha, \beta, \gamma)] \in \tilde{\Sigma}^3}
(a_{\alpha} a_{\beta} + a_{\beta} a_{\gamma} +
a_{\gamma} a_{\alpha}) N_{\alpha, \beta} \ea \wedge E_{\beta} \wedge
E_{\gamma}
\end{eqnarray}
Clearly, $ih_{\alpha} \wedge \ea \wedge \eb \in \wedge^3 \fk$ for each $\a
\in \Sigma_{+}$. 
Suppose that $(\alpha, \beta, \gamma) \in \Sigma^3$ are such that
$\alpha + \beta + \gamma = 0$. Without loss of generality, we
can assume that $\a, \beta \in \Sigma_{+}$ and $\gamma  \in -\Sigma_{+}$.
Then
\[
N_{\alpha, \beta} \ea \wedge E_{\beta} \wedge
E_{\gamma} \, + \, N_{-\alpha, -\beta} \eb \wedge E_{-\beta} \wedge
E_{-\gamma} \, = \, N_{\alpha, \beta} (\ea \wedge E_{\beta} \wedge
E_{\gamma} - \eb \wedge E_{-\beta} \wedge E_{-\gamma} ).
\]
This element is in $\wedge^3 \fk$ because it is fixed by $\theta \in
{\rm End}_{{\Bbb R}}(\wedge^3 \fg)$ defined by
\[
\theta (x_1 \wedge x_2 \wedge x_3) \, = \, \theta(x_1) 
\wedge \theta(x_2) \wedge
\theta (x_3), \hspace{.3in} x_1, x_2, x_3 \in \fg,
\]
where $\theta 
\in {\rm End}_{{\Bbb R}} (\fg)$ is the complex conjugation of $\fg$
defined by $\fk$. The right hand side of (\ref{eq_in-k}) is thus
the Schouten bracket of $A$ with itself inside $\wedge \fk$.
 
Now suppose that $r$ is a classical dynamical $r$-matrix
for the pair $(\fg, \fh)$ as given in Theorem \ref{thm_ev}.
Suppose that $\lambda \in \fh^*$ is in the domain of $r$ such that
the skew-symmetric part $A_r(\lambda) = r(\lambda) -
{\epe \over 2} \Omega$ of $r(\lambda)$ lies in $\wedge^2 \fk$.
Then
\[
[A_r(\lambda), \, A_r(\lambda)] \, - \, [\Lambda, \, \Lambda]
\in (\wedge^3 \fk) \cap (\fh \wedge \fk \wedge \fk) =
\ft \wedge \fk \wedge \fk.
\]
By abuse of notation, we still use $\tilde{\pi}_{r(\lambda)}$
(already used in Theorem \ref{thm_main}) to denote the bi-vector field
on $K$ given by
\[
\tilde{\pi}_{r(\lambda)} (k) \, = \, R_k \Lambda \, - \, L_k A_r(\lambda),
\hspace{.2in} k \in K,
\]
where $R_k$ and $L_k$ are respectively the right and left
translations on $K$ by $k$. We use $\pi_{r(\lambda)}$
to denote the projection of $\tilde{\pi}_{r(\lambda)}$ to $K/T$ by the map
$p: K \rightarrow  K/T: k \mapsto kT$.
 
\begin{thm}
\label{thm_compact}
 Let $r$ be any classical dynamical $r$-matrix for
the pair $(\fg, \fh)$ given in Theorem \ref{thm_ev}.
Suppose that $\lambda \in \fh^*$ is in the domain of $r$ such that
$A_r(\lambda) = r(\lambda) - {\epe \over 2} \Omega$
is in $\wedge^2 \fk$. Then,

1) the bi-vector field
$\pi_{r(\lambda)}$ on $K/T$ defines a $(K, \piK)$-homogeneous
Poisson structure on $K/T$;
 
2) with the Poisson structure $\piK$  on $K$ given by (\ref{eq_on-K}),
every $(K, \piK)$-homogeneous Poisson structure on $K/T$
arises this way.
\end{thm}
 
\noindent
{\bf Proof.} The proof of 1) is similar to that of Theorem \ref{thm_main}.
We prove 2). Assume that $\pi$ is a $(K, \piK)$-homogeneous
Poisson structure
on $K/T$. Since $\pi$ is $T$-invariant, we can write
\[
\pi(e) \, = \, \sum_{\alpha \in \Sigma_{+}} 
(-{i \epe \over 2} + i \phi_{\alpha})
{\Xa \wedge \Ya \over 2} \in \wedge^2 (\fk / \ft),
\]
where $e = eT \in K/T$ and $\phi_{\alpha} \in i {\Bbb R}$ for each $\alpha
\in \Sigma_{+}$. 
(Recall that $\epe \in i {\Bbb R}$ is fixed at the beginning.)
Set $\phi_{-\alpha } = - \phi_{\alpha}$ for $\a \in \Sigma_{+}$. 
Using the same trick
for calculating the Schouten bracket in $\wedge \fk$, i.e., by
embedding $\wedge \fk$ into $\wedge \fg$, and by using arguments similar to
those in the proof of Lemma \ref{lem_three}, we know that the
$\phi_{\alpha}$'s must satisfies Condition (\ref{eq_phi}). Exactly the
same as in the proof of the second part of Theorem \ref{thm_main},
we know that there exist a choice of positive roots $\Sigma^{'}_{+}$, 
a choice
of a subset $X^{'}$ of the set of simple roots for $\Sigma^{'}_{+}$, and
some (not necessarily unique) $\lambda_0 \in \fh^*$ such that
\[
\phi_{\alpha} \, = \, \left\{ \begin{array}{ll}
{\frac{\epe}{2}} \coth {\frac{\epe}{2}} \ll \alpha, \, \lambda_0
\gg & {\rm if} \alpha \in [X^{'}] \\
\pm {\epe \over 2} & {\rm if} \alpha \in \pm (\Sigma^{'}_{+} \backslash
[X^{'}].
\end{array} \right.
\]
Let $r$ be the classical dynamical $r$-matrix for the pair $(\fg, \fh)$
defined by $\Sigma^{'}_{+}$ and $X^{'}$ 
as in Theorem \ref{thm_ev} ($\mu = 0$ and
$C = 0$), we see that $\pi$ coincides with the Poisson structure
$\pi_{r(\lambda_0)}$ on $K/T$.
\qed

\section{The Poisson structures $\pix$ on $K/T$}
\label{sec_rX}

\subsection{Definition}
\label{sec_dfn-pix}

As in the case for $G/H$, we will single out a family of 
$(K, \piK)$-homogeneous Poisson structures
on $K/T$ which exhausts all such Poisson structures on $K/T$
up to $K$-equivariant isomorphisms.

For a subset $X \subset S(\Sigma_+)$,  set 
\[
\fa_{\tx} \, = \, {\rm span}_{{\Bbb R}} \{h_{\gamma} = [E_\gamma, 
E_{-\gamma}]: \gamma \in X\}.
\]
Denote by $\{\check{h}_{\gamma}: \gamma \in S(\Sigma_{+}) \}$
the set of fundamental co-weights for  $S(\Sigma_{+})$, i.e.,
$\check{h}_{\gamma} \in \fa$ for each $\gamma \in S(\Sigma_+)$
and $\gamma_1(\check{h}_{\gamma}) = \delta_{\gamma_1, \gamma}$
for all $\gamma_1, \gamma \in S(\Sigma_+)$..
For $X_1 \subset S(\Sigma_+)$, set
\[
\check{\rho}_{\txo} \,  =\,  \sum_{\gamma \in \txo} \check{h}_{\gamma}.
\]
Define $\check{\rho}_{\txo}$ to be  $0$ if
$X_1$ is the empty set.

\begin{thm}
\label{thm_pix}
For $X \in S(\Sigma_{+}), X_1 \subset X$ and $\lambda 
= \lambda_1 + {i \pi \over 2} \check{\rho}_{\txo}
\in 
\fa_{\tx} + {i \pi \over 2} \check{\rho}_{\txo}$
such that $\alpha (\lambda_1) \neq 0 $ for
all $\alpha \in [X]$ with $\alpha(\check{\rho}_{\txo})$ even, let
$\pix$ be the bi-vector field on $K/T$ given by
\[
\pix \, = \, p_{*} \piK \, - \, {i \epe \over 2} \left(
\sum_{\alpha \in [X] \cap \Sigma_+} {1 \over 1 - e^{2 \alpha(\lambda)}}
\Xa \wedge \Ya \right)^L,
\]
where the second term 
on the right hand side is the $K$-invariant bi-vector field on $K/T$
whose value at $e = eT$ is the expression given in the parenthesis.
Then

1) $\pix$ is a $(K, \piK)$-homogeneous Poisson structure
on $K/T$, and

2) every $(K, \piK)$-homogeneous Poisson 
structure on $K/T$ is $K$-equivariantly isomorphic 
to some $\pix$.
\end{thm}

\begin{rem}
\label{rem_equi}
{\em
Note that the condition on $\lambda_1 \in \fa_{\tx}$ is equivalent to 
$\alpha (\lambda) \notin \pi i \Z$ for all $\alpha \in [X]$, so that
$e^{2\alpha(\lambda)} \neq 1$ for all $\alpha \in [X]$.
}
\end{rem}

\noindent
{\bf Proof.}
1). The number $e^{2\alpha(\lambda)}$ is real for each $\alpha \in [X]$. 
Thus $\pix$ is a $(K, \piK)$-homogeneous Poisson structure
coming from a classical dynamical $r$-matrix.

2) Assume that $\pi$ is a $(K, \piK)$-homogeneous Poisson 
structure on $K/T$. By Theorem \ref{thm_compact} and by a proof similar to
that of Theorem \ref{thm_G}, there exist
$X \subset S(\Sigma_{+})$ and some $\lambda_0 \in \fh^*$ such that
$\pi$ is isomorphic, via a $K$-equivariant diffeomorphism of $K/T$,
to the Poisson structure $\pi^{'}$ given by
\[
\pi^{'} \, = \, p_{*} \piK \, - \, {i \epe \over 2} 
\left(\sum_{\alpha \in [X] \cap \Sigma_{+}} 
k_{\alpha} \Xa \wedge \Ya\right)^L,
\]
where
\[
k_{\alpha} \, = \, {\frac{1}{2}} (1 - \coth ({\frac{\epe}{2}} \ll
\alpha, \lambda_0 \gg)) \, = \, 
{1 \over 1 - e^{\varepsilon \ll \alpha, \lambda_0 \gg}} \in \R.
\]
Let $h_{\lambda_0} \in \fh$ be the element in $\fh$ corresponding to 
$ \lambda_0$ under the Killing form, so that 
$\ll \alpha, \lambda_0 \gg = \alpha(h_{\lambda_0})$
for all $\alpha \in \Sigma$.
It remains to show that ${\epe \over 2} h_{\lambda_0}$
can be replaced by some $\lambda  \in
\fa_{\tx} + {i \pi \over 2} \check{\rho}_{\txo}$. To this end,
consider the function $f(z) = 1/(1-e^z)$ for $z \in \C$. It takes
values in all of $\C$ except for $0$ and $1$. Moreover, 
$f(\R \backslash \{0\}) = (-\infty, 0) \cup (1, \infty)$ and
$f(\R + i \pi) \in (0,1)$. Set
\[
X_1 \, = \, \{\gamma \in X: \, \, k_{\gamma} \in (0,1)\}.
\]
Then for each  $\gamma \in X$, there exists $\mu_{\gamma}
\in \R$ such that
\[
\left\{ \begin{array}{ll} & k_{\gamma} \, = \, 
f(\mu_{\gamma} + i \pi) \hspace{.2in} {\rm if} \hspace{.1in}
\gamma \in X_1 \\
& k_{\gamma} \, = \,
f(\mu_{\gamma} ) \hspace{.2in} {\rm if} \hspace{.1in} \gamma \in 
X \backslash X_1. \end{array} \right.
\]
Let $\lambda_1 \in \fa_{\tx}$ be such
that $2\gamma(\lambda_1) = \mu_\gamma$ for each 
$\gamma \in X$ (such a $\lambda_1$ exists), and
let
$\lambda = \lambda_1 + {\pi i \over 2}
\check{\rho}_{\txo}$.
Then
$k_{\gamma} = f(2\gamma(\lambda))$ for all $\gamma \in X$.
Consequently, by writing $\alpha \in [X] \cap \Sigma_+$ as a linear
combination of elements in $X$, we see that
 $k_{\alpha} = f(2\alpha((\lambda))$ for all
$\alpha \in [X]$.
\qed

\begin{nota}
\label{nota_infty}
{\em For reasons given in Section \ref{sec_limits}, we will 
use $\pinf$ to denote the Poisson structure
$p_* \piK$ on $K/T$. It is called the 
{\it Bruhat Poisson structure} \cite{lu-we:poi}
because its symplectic leaves are Bruhat cells in $K/T$.
See Section \ref{sec_leaves-1} for more details.
}
\end{nota}

\begin{exam}
\label{exam_sl2}
{\em
Consider 
\[
K \, =\,  SU(2) \, = \, \left\{
\left( \begin{array}{ll} u & v \\ -\bar{v} &  \bar{u} 
\end{array} \right): \, \, u, v \in \C, \, |u|^2 + |v|^2 = 1\right\},
\]
$T = \{{\rm diag}(
e^{ix}, e^{-ix}): \, x \in \R\} \cong S^1$ and the root 
$\alpha(x, -x) = 2x$ is taken to be the positive root. Then
\[
\Xa \, = \, {\frac{1}{2}} \left( \begin{array}{ll} 0 & 1 \\ -1 & 0
\end{array} \right), \hspace{.2in}
\Ya \, = \, {\frac{1}{2}} \left( \begin{array}{ll} 0 & i \\  i & 0
\end{array} \right).
\]
With 
\[
\Lambda \, = \, -{i \epe \over 2} {\Xa \wedge \Ya \over 2} \, \in \, 
\su(2) \wedge \su(2)
\]
and the Poisson structure $\piK$ on $K = SU(2)$ defined by
\[
\piK \, = \, \Lambda^R \, - \, \Lambda^L,
\]
the Poisson brackets among the coordinate functions $u, v, \bar{u}$
and $\bar{v}$ on $SU(2)$ are given by
\[
\{u, \, \bar{u}\} = -{\epe \over 4} |v|^2, \hspace{.2in}
\{u, \, v\} = {\epe \over 8} uv, \hspace{.2in}
\{u, \, \bar{v}\} = {\epe \over 8} u \bar{v}, \hspace{.2in}
\{v, \bar{v}\} = 0.
\]
Let $\pi_0$ be the $SU(2)$-invariant bivector
field on $SU(2)/S^1$ whose value at the point $e = eS^1$ is
$\Xa \wedge \Ya$.
It is symplectic.

{\bf Case 1}: $X = X_1 = \emptyset$. Then $\pix = \pinf$;

{\bf Case 2}: $X = \{\alpha\}, \, X_1 = \emptyset$. Then
$\lambda = \left(\begin{array}{cc} \lambda_1 & 0 \\ 0 & -\lambda_1\end{array}
\right)$ with $\lambda_1 \neq 0$, and
\[
\pix \, = \, \pinf \, - \, {i \epe \over 2}
 {1 \over 1 - e^{4 \lambda_1}} \pi_0.
\]
 
{\bf Case 3:} $X = X_1 = \{\alpha\}$. Then
\[
\lambda \, = \,
\left(\begin{array}{cc} \lambda_1 + {\pi i \over 4}
 & 0 \\ 0 & -\lambda_1 - {\pi i \over 4}
\end{array}\right)
\]
with $\lambda_1 \in \R$ arbitrary, and
\[
\pix \, = \, \pinf \, - \, {i\epe \over 2}
 {1 \over 1 + e^{4 \lambda_1}} \pi_0.
\]
Note that the range of the function ${1 \over 1 - e^{4 \lambda_1}}$
for $\lambda_1 \in \R\backslash\{0\}$ is $(-\infty, 0) \cup (1, +\infty)$,
and the range of ${1 \over 1 + e^{4 \lambda_1}}$ for $\lambda_1 \in \R$ is
$(0, 1)$. Thus, for all possible choices of $X, X_1$ and $\lambda$,
we get all the Poisson structures of the form
\[
\pi^a \, = \, \pinf \, - \, {i \epe \over 2} a \pi_0
\]
for $a \in \R$
except for $a = 1$.
But the Poisson structure $\pi^a$
when $a = 1$ is easily seen to be isomorphic to $\pinf$ (corresponding
to $a = 0$)
by the $SU(2)$-equivariant diffeomorphism on $SU(2)/S^1$
defined by the right translation by the non-trivial Weyl group
element. 
The fact that every $(SU(2), \pi_{\tk})$-homogeneous
Poisson structures on $S^2$ is of the form $\pi^a$ for some $a \in \R$
is very easy to  check directly \cite{sheu:s2}.
 
Identify  the Lie algebra $\su(2)$ with $\R^3$ by
\[
\left( \begin{array}{ll} ix & y + iz\\-y+iz & -ix\end{array} \right)
\, \Map \, (x, y, z)
\]
so the Adjoint orbit through $\left( \begin{array}{ll} i & 0 \\ 0 & -i
\end{array} \right)$ can be identified with the sphere $S^2
= \{(x, y, z) \in \R^3: x^2 + y^2 + z^2 = 1\}$.
Consequently, we have the identification
\[
SU(2) /S^1 \rightarrow S^2: \, kS^1 \Map {\rm Ad}_k
\left( \begin{array}{ll} i & 0 \\ 0 & -i
\end{array} \right),
\]
or
\[
\left( \begin{array}{ll} u & v \\ -\bar{v} &  \bar{u}
\end{array} \right) S^1 \Map (|u|^2 - |v|^2, \, -i(uv - \bar{u} \bar{v}),
\, -(uv + \bar{u} \bar{v}) ).
\]
The induced Bruhat Poisson structure $\pinf$ on $S^2$
is given by
\[
\{x, y\} = -{\epe i \over 4} (x-1) z, \hspace{.2in}
\{y, z\} = -{\epe i \over 4} (x-1) x, \hspace{.2in}
\{z, x\} = -{\epe i \over 4} (x-1)y,
\]
and the Poisson structure $\pi^a$ on $S^2$ is given by
\[
\{x, y\} = -{\epe i \over 4} (x+2a-1) z, \hspace{.2in}
\{y, z\} = -{\epe i \over 4} (x+2a-1) x, \hspace{.2in}
\{z, x\} = -{\epe i \over 4} (x+2a-1)y,
\]
Note that $\pi^a$ is symplectic when $a < 0$ or $a > 1$. When
$a = 0$, it has two symplectic leaves, the point $(1, 0, 0)$
being a one-point leaf and the rest of $S^2$ as another leaf.
Similarly for $a = 1$. When $0 < a < 1$, it has infinitely many symplectic
leaves: two open leaves respectively given by $x < 1-2a$ and
$x > 1-2a$, and every point on the circle $x = 1-2a$ as
a one-point leaf.
}
\end{exam}

\begin{exam}
\label{exam_su3}
{\em
Let $\fg = \fsl(3, \C)$ and $K = SU(3)$. The
three positive roots are chosen to be
\[
\alpha_1 (x) = x_1 - x_2, \hspace{.2in}
\alpha_2 (x)= x_2 - x_3, \hspace{.2in}
\alpha_3(x) = x_1 - x_3
\]
for a diagonal matrix $x = {\rm diag}(x_1, x_2, x_3)$. 
Take $X = S(\Sigma_{+}) = \{\alpha_1, \alpha_2\}$ and
$X_1 =  \{ \alpha_1\}$. In this case
\[
\check{\rho}_{\txo} \, = \, \left(\begin{array}{ccc} {2 \over 3} & 0 & 0 \\
0 & -{1 \over 3} & 0\\
0 & 0 & -{1 \over 3} \end{array}\right),
\]
and
\[
\lambda \, = \, \left(\begin{array}{ccc} \lambda_1 + {\pi i \over 3} & 0 & \\
0 & \lambda_2 - {\pi i \over 6} & 0 \\
0 & 0 & -(\lambda_1 + \lambda_2) - {\pi i \over 6} \end{array}
\right), \hspace{.1in} \lambda_1 + 2\lambda_2 \neq 0.
\]
Then
\[
\pix  =   \pinf  +\left(
{2 X_{\alpha_1} \wedge Y_{\alpha_1} \over 1 + e^{2(\lambda_1 - \lambda_2)}} 
 \,+\,
{2X_{\alpha_2} \wedge Y_{\alpha_2} \over 1 - e^{2\lambda_1 +4 \lambda_2}} 
 \,+\,
{2 X_{\alpha_3} \wedge Y_{\alpha_3} \over 1 + e^{4\lambda_1 +2 \lambda_2}} 
\right)^L.
\]
 }
\end{exam}

\subsection{Connections via taking limits in $\lambda$}
\label{sec_limits}

As noted in \cite{e-v:cdyb}, the dynamical $r$-matrices
are related to each other via taking various limits in $\lambda$.
Correspondingly, the Poisson structures $\pix$ are also related
this way. We study these relations in the section.

\begin{prop}
\label{prop_limit}
For any $X_1 \subset X \subset Y \subset S(\Sigma_{+})$ and
$\lambda = \lambda_1 + {i \pi \over 2} \check{\rho}_{\txo}
 \in \fa_{\tx} + 
{i \pi \over 2} \check{\rho}_{\txo}$ such that
$\alpha(\lambda_1) \neq 0$ for all $\alpha \in [X]$ with
$\alpha(\check{\rho}_{\txo})$ even,
we have
\begin{equation}
\label{eq_limit}
\pix \, = \, \lim_{t \rightarrow +\infty} \pi_{\ty, \txo, \lambda+ t 
\check{\rho}_{\ty \backslash \tx}}.
\end{equation}
In particular, 
\[
\pinf \, = \, \lim_{t \rightarrow +\infty} \pi_{\ty, \emptyset, 
t \check{\rho}_{\ty}}.
\]
Moreover, we also have
\begin{equation}
\label{eq_pinf}
\pinf \, = \, \lim_{t \rightarrow +\infty} \pi_{
\tx, \txo, \lambda + t \check{\rho}_{\tx}}.
\end{equation}
\end{prop}

\noindent
{\bf Proof.} Set $\mu_t = \lambda+ t\check{\rho}_{\ty \backslash \tx}$
for $t > 0$. Let $\a \in [Y] \cap \Sigma_+$. If $\a \in [X]$, then
$\a(\check{\rho}_{\ty \backslash \tx}) = 0$ so
$\a (\mu_t) = \a(\lambda)$. If
$\a \in [Y] \backslash [X]$, then $v := 
\a(\check{\rho}_{\ty \backslash \tx})$ is positive, so
\[
\lim_{t \rightarrow \infty} {1 \over 1 - e^{\alpha(\mu_t)}}
\, = \, \lim_{t \rightarrow \infty} {1 \over 1 - e^{tv}} \, = \, 0.
\]
Hence (\ref{eq_limit}) follows from the definition of $\pix$.
The limit in (\ref{eq_pinf}) is obvious.
\qed

\subsection{The Lagrangian subalgebras of $\fg$ corresponding to
$\pix$}
\label{sec_lagrangian-compact}

The Lie bialgebra of the Poisson Lie group $(K, \piK)$
is $(\fk, \fa + \fn)$, where the pairing between
$\fk$ and $\fa + \fn$ is given by 
${\frac{2i}{\varepsilon}} {\rm Im} \ll \, , \, \gg$, where
${\rm Im} \ll \, , \, \gg$ stands for the imaginary part of
the Killing form $\ll \, , \, \gg$.

We will call a real subalgebra $\fl$ of $\fg$ a {\it Lagrangian algebra} 
if
1) $\dim \fl = \dim \fk$, and 2) ${\frac{2i}{\varepsilon}} {\rm Im} \ll
x, y\gg = 0$ for all $x, y \in \fl$. By a theorem
of Drinfeld \cite{dr:homog},  $(K, \piK)$-homogeneous
Poisson structures on $K/T$ correspond to Lagrangian 
subalgebras $\fl$ of $\fg$ with $\fl \cap \fk = \ft$. 
In this section, we calculate the Lagrangian subalgebras 
$\lix$ corresponding
to the Poisson structures $\pix$.

By definition \cite{dr:homog}, 
\[
\lix \, = \, \{x + \xi: \, x \in \fk, \, \xi \in 
\fa + \fn: \, \xi|_{\frak t} = 0, \, 
\xi \backl \pix(e) = x + \ft\}.
\]
A direct calculation gives
\beqa
\lix & = & \ft \, + \, {\rm span}_{\Bbb R} \{
E_{\beta}, iE_{\beta}: \, \beta \in \Sigma_{+} \backslash [X]\}\\
& & + {\rm span}_{\Bbb R} \{ {1 \over e^{2\alpha(\lambda)} -1} \Xa + \ea, 
\, \, \, 
{1 \over e^{2\alpha(\lambda)} -1} \Ya + i\ea: \, \,\alpha \in [X]
\cap \Sigma_+\}.
\eeqa
On the other hand, for $\alpha \in [X]$, since $e^{2\alpha(\lambda)}
\neq 1$, we have
\beqa
{\rm Ad}_{e^{\lambda}} \Xa & = & {\rm Ad}_{e^{\lambda}} (\ea - \eb)
\, = \, (e^{\alpha (\lambda)} - e^{-\alpha (\lambda)} )
({1 \over e^{2 \alpha (\lambda)} -1} \Xa + \ea) \\
{\rm Ad}_{e^{\lambda}} \Ya & = & {\rm Ad}_{e^{\lambda}} (i\ea +i \eb)
\, = \, (e^{\alpha (\lambda)} - e^{-\alpha (\lambda)} )
({1 \over e^{2 \alpha (\lambda)} -1} \Ya + i\ea).
\eeqa
Note that $e^{\alpha(\lambda)}$ is real or imaginary depending 
on $\alpha(\check{\rho}_{\txo})$ is even or odd.
Set
\begin{equation}
\label{eq_mx}
\fn_{\tx} \, = \, {\rm span}_{\Bbb R} \{
E_{\beta}, iE_{\beta}: \beta \in \Sigma_{+} \backslash [X]\}.
\end{equation}
Then we have proved
the following proposition.

\begin{prop}
\label{prop_lix}
Denote by $\lix$ the Lagrangian subalgebra of $\fg$ 
corresponding to the Poisson
structure $\pix$ on $K/T$. It is given by
\beqa
\lix  \, =\,  {\rm Ad}_{e^{\lambda}} (
\ft\, + \, \fn_{\tx} &   + &
{\em span}_{\Bbb R} \{\Xa, \Ya: \, \alpha \in [X], \,  
\alpha(\check{\rho}_{\txo}) \, {\em is} \,\, {\em even}\}\\
 & + &  {\em span}_{\Bbb R} \{i\Xa, i\Ya: \, 
\alpha \in [X], \, \alpha(\check{\rho}_{\txo}) \, 
{\em is} \,\, {\em  odd}\} ).
\eeqa
\end{prop}

\begin{rem}
\label{rem_signature}
{\em
Let $\theta$ be the complex conjugation on $\fg$ defined by 
$\fk$. Let $\tau_{\tx, \txo}$ be the complex conjugation on $\fg$
given by
\[
\tau_{\tx, \txo} \, = \, {\rm Ad}_{\exp(\pi i \check{\rho}_{\txo})}
\theta
\, = \, \theta {\rm Ad}_{\exp(-\pi i \check{\rho}_{\txo})}.
\]
Denote by $\fm_{\tx}^{\tau_{\tx, \txo}}$ the set of fixed
points of 
$\tau_{\tx, \txo}$ in $\fm_{\tx}$, where
\[
\fm_{\tx} \, = \, \fh \, + \, 
{\rm span}_{\Bbb C} \{\ea: \alpha \in [X]\}.
\]
Then
\[
\lix \, = \, {\rm Ad}_{e^{\lambda}} (\fm_{\tx}^{\tau_{\tx, \txo}} + 
\fn_{\tx}).
\]
}
\end{rem}

\begin{rem}
\label{rem_grass}
{\em
Let $n = \dim \fk$ and consider $\lix$ as a point in
$\Gr(n, \fg)$, the Grassmannian of $n$-dimensional real subspaces
of $\fg$. Then, corresponding to Proposition \ref{prop_limit},
we have, for $X_1 \subset X \subset Y \subset S(\Sigma_{+})$ 
and for any $\lambda = \lambda_1 + {i \pi \over 2}
\check{\rho}_{\txo} \in \fa_{\tx} + {i \pi \over 2}
\check{\rho}_{\txo}$ such that
$\alpha(\lambda_1) \neq 0$ for all $\alpha \in [X]$ with 
$\alpha(\check{\rho}_{\txo})$ even,
\begin{equation}
\label{eq_grass}
\lim_{t \rightarrow + \infty} \fl_{\ty, \txo, \lambda +
t \check{\rho}_{\ty \backslash \tx}} \, = \, \lix
\end{equation}
in $\Gr(n, \fg)$. Indeed, under the Plucker embedding of
$\Gr(n, \fg)$ into ${\Bbb P}^1 (\wedge^n \fg)$, the
Lie subalgebra $\liy$
corresponds to the point in ${\Bbb P}^1 (\wedge^n \fg)$ defined by the vector
\[
v_{\ty, \txo, \lambda} := Z_0 \wedge \prod_{\alpha \in [\ty] \cap \Sigma_+}
\left( {1 \over e^{2 \alpha (\lambda)} -1} \Xa + \ea \right)
\wedge \left({1 \over e^{2 \alpha (\lambda)} -1} \Ya + i\ea \right)
\wedge \prod_{\alpha \in \Sigma_{+} \backslash [\ty]} \ea \wedge \eb
\]
where $Z_0 \in \wedge^{\dim \ft} \ft$ and $Z_0 \neq 0$ is fixed.
Since $v_{\ty, \txo, \lambda + t\check{\rho}_{\ty \backslash \tx}}
\rightarrow v_{\txo, \lambda}$ as $t \rightarrow +\infty$, we
see that (\ref{eq_grass}) holds in ${\Bbb P}^1 (\wedge^n \fg)$
and thus also in $\Gr(n, \fg)$. 
}
\end{rem}
 
\begin{exam}
\label{exam_lix-X1-empty}
{\em
When $X = X_{1}$ are the empty set, we have 
$\lix = \ft + \fn$, and
when $X = S(\Sigma_{+})$ and $X_{1}$ is the empty set, we have
$\lix = {\rm Ad}_{e^{\lambda}} \fk$. In general, 
when $X = S(\Sigma_{+})$, the Lie subalgebra $\lix$ is a real
form of $\fg$.
}
\end{exam}

\subsection{Geometrical interpretation of $\pix$}
\label{sec_geom-X-whole}

Denote by $\Lagr$ the set of all Lagrangian subalgebras of $\fg$
with respect to the imaginary part of the Killing form
$\ll \, , \, \gg$. (Here $\fg$ is regarded as a real
vector space.)
It is an algebraic subvariety of the Grassmannian 
$\Gr(n, \fg)$ of $n$-dimensional subspaces of $\fg$, where $n = \dim \fk$.
In \cite{e-l:Lagrangian}, we show that there is a smooth bivector
field $\Pi$ on $Gr(n, \fg)$ such that the Schouten bracket
$[\Pi, \Pi]$ vanishes at every $\fl \in \Lagr$.
More precisely,
consider the $G$-action on $\Gr(n, \fg)$ by the Adjoint action. 
It defines a Lie algebra anti-homomorphism 
\[
\kappa: \, \fg \lrw \chi^1(\Gr(n, \fg)),
\]
where $\chi^1(\Gr(n, \fg))$ is the space of vector fields on 
$\Gr(n, \fg)$. Denote by the same letter its 
multi-linear extension from $\wedge^2 \fg$ to the 
space of bi-vector fields on $\Gr(n, \fg)$. Then the
bivector field $\Pi$ on $\Gr(n, \fg)$ is defined to be
\[
\Pi \, = \, {1 \over 2} \kappa(R),
\]
where $R \in \wedge^2 \fg$ is the $r$-matrix for $\fg$ given by
\begin{equation}
\label{eq_R}
\la\, R, \,\, (x_1 + y_1) \wedge (x_2 + y_2)\, \ra_{\varepsilon}
\, = \,
\la x_1, \, y_2 \ra_{\varepsilon} \, - \, \la x_2, \, y_1
\ra_{\varepsilon}
\end{equation}
for $x_1, x_2 \in \fk$  and $y_1, y_2 \in \fa + \fn$ with
$\la \, , \, \ra_{\varepsilon} = {2i \over \epe} {\rm Im} \ll \, , \, 
\gg$. Explicitly,
\[
R \, = \,  - {\epe \over 2i} \left(
\sum_{j=1}^{l} (i h_j) \wedge h_j \, + \, \sum_{\alpha \in \Sigma_{+}}
(-\Xa \wedge (i \ea) + \Ya \wedge \ea) \right),
\]
where $\{h_1, ..., h_l\}$ is a basis for $\fa$ such that
$\ll h_j, h_k \gg = \delta_{jk}$. 
It now follows from the definition of $\Pi$ that
it defines a Poisson structure on every $G$-invariant
smooth submanifold of $\Lagr$.

One particular $G$-invariant smooth submanifold of $\Lagr$
is the (unique) irreducible component $\Lagr_0$ of $\Lagr$
that contains $\fk$. We show in \cite{e-l:Lagrangian} that
each $\lix \in \Lagr_0$ and that its $K$-orbit in $\Lagr_0$
is a Poisson submanifold of $(\Lagr_0, \Pi)$.
(We also show in \cite{e-l:Lagrangian} that
$\Lagr_0$ is diffeomorphic to the set of real points in
the De Concini-Procesi compactification of $G$ 
\cite{dp:compactification}). For each Poisson structure
$\pix$ on $K/T$, consider the map
\[
P: \,  (K/T, \, \pix) \lrw (\Lagr_0, \, \Pi): \, 
kT \, \Map \, {\rm Ad}_k \lix.
\]
It is shown in \cite{e-l:Lagrangian} that $P$ is
a Poisson map. When the normalizer subgroup of $\lix$
in $K$ is $T$, this map is an embedding of 
$K/T$ into $\Lagr_0$ whose image is the 
the $K$-orbit of $\lix$ in $\Lagr_0$. In general,
$P$ is a covering map onto the $K$-orbit of $\lix$ in $\Lagr_0$.
Thus, every $(K/T, \pix)$ is a Poisson
submanifold of $(\Lagr_0, \Pi)$ (possibly up to a covering
map).
This can be considered as one geometrical interpretation
of $\pix$.

Two special cases of $\pix$ deserve more attention. The
first is when $X = X_1 = \emptyset$ ($\lambda = 0$ in this case).
Then $\pix = \pinf$ is the Bruhat Poisson structure. It has been the 
most interesting example in terms of connections to Lie theory.
For its relations with the Kostant harmonic forms \cite{ko:63}, see
\cite{lu:coor} and \cite{e-l:harm}.

The second special case is when 
$X = S(\Sigma_+)$
and $X_1 = \emptyset$. The condition on $\lambda
$ is that $\lambda \in \fa$ is regular.  We will show that
$\pix$ is symplectic in this case. In fact, we will show that 
$\pix$ can be identified with the symplectic structure
on a dressing orbit of $K$ in its dual Poisson Lie group. 
We also remark that this symplectic structure
has been used in \cite{lu-ra:convexity} to give a symplectic
proof of Kostant's nonlinear convexity theorem.

Recall that 
the Manin triple
$(\fg, \fk, \fa + \fn,  {2i \over \epe} {\rm Im} \ll\, , \, \gg)$
gives rise to a Poisson structure $\pi_{AN}$ on the group $AN$
making $(AN, \pi_{AN})$ into the dual Poisson Lie group
of $(K, \piK)$.
The
group $K$ acts on $AN$ by the {\it (left) dressing action}:
\[
K \times AN \lrw AN: \, \, (k, b) \Map k \cdot b :=b_1, \hspace{.3in}
{\rm if} \, \, b k^{-1} = k_1 b_1 \,\, {\rm for}\, \, k_1 \in K
\,\, {\rm and}\, \, b_1 \in AN.
\]
The $K$ orbits of this dressing action of $K$ in $AN$, called the
{\it dressing orbits}, are precisely all the symplectic leaves
of the Poisson structure
on $AN$ and they are parametrized by a fundamental $W$-chamber in $\fa$.
Thus each dressing orbit inherits a symplectic, and thus Poisson,
structure as a symplectic leaf. Since the dressing action
is Poisson \cite{sts:dressing} \cite{lu-we:poi}, these
dressing orbits are examples of
$(K,  \piK)$-homogeneous Poisson spaces.
Let $\lambda \in \fa$ be regular and consider the element
$\el \in A$.
The stabilizer subgroup of $K$ in $AN$
at $\el$ is $T$. Thus, by identifying $K/T$ with the dressing orbit
through $\el$, we get a Poisson structure on $K/T$ which is in fact
symplectic.

\begin{nota}
\label{nota_symplectic}
{\em
We will use
$\pi_{\lambda}$ to denote
the Poisson structure on $K/T$ obtained by identifying
$K/T$ with the symplectic leaf in $AN$ through the point $\el$, and
we call it the
{\it dressing orbit Poisson structure corresponding to $\el \in A$}.
}
\end{nota}

\begin{prop}
\label{prop_dressing}
When $X = S(\Sigma_{+}), X_1 = \emptyset,
$ and $\lambda \in \fa$ is regular, the
Poisson structure $\pix$ on $K/T$ is nothing but
the dressing orbit Poisson structure $\pi_{\lambda}$ corresponding to
$\el$. Explicitly, we have
\begin{equation}
\label{eq_pi-lambda}
\pi_{\lambda} \, = \, -{\frac{i \epe }{2}}
\left( \sum_{\alpha \in \Sigma_{+}}
{\frac{1}{1-e^{2\alpha(\lambda)}}}
\Xa \wedge \Ya \right)^{L} \, + \, \pi_{\infty},
\end{equation}
where the first term is the $K$-invariant bi-vector field on $K/T$
whose value at $e = eT$ is the expression given in the parenthesis.
\end{prop}

\noindent
{\bf Proof.} Since $\lix$ is given by the right hand side of 
(\ref{eq_pi-lambda}), we
only need
to show that the dressing orbit Poisson structure $\pi_{\lambda}$
is also given by the same formula.
Denote the Poisson structure on $AN$ by $\pi_{\scriptscriptstyle AN}$.
Since we are identifying $\fk$ with $(\fa + \fn)^*$ via
${\frac{2i}{\epsilon}} {\rm Im} \ll \, , \, \gg$, an element $x \in \fk$
can be
regarded as a left invariant $1$-form on $AN$ which we denote by
$x^l$. Let $p_{\frak k}: \fg \rightarrow \fk$ be
the projection from $\fg$ to $\fk$
with respect to the Iwasawa Decomposition $\fg = \fk + \fa + \fn$.
We know that (see \cite{lu:thesis}) for any $a \in A$,
\[
\pi_{\scriptscriptstyle AN}(x^l, y^l) (a) \, = \, {2i \over \epe} {\rm Im}
\ll {\rm Ad}_a x, \, p_{\frak k} {\rm Ad}_a y \gg
\]
for all $x, y \in \fk$. Here, ${\rm Ad}_a$ is the Adjoint action of
$a \in A$ on $\fg$. Thus, when $x$ and $y$ run over the basis
vectors $\{iH_{\alpha}, \Xa, \Ya: \a \in \Sigma_{+} \}$ for $\fk$,
we have $\pi_{\scriptscriptstyle AN}(x^l, y^l) (a) = 0$ except that
\beqa
\pi_{\scriptscriptstyle AN} (X_{\alpha}^{l}, Y_{\alpha}^{l}) & = &
{\frac{2i}{ \epsilon}} {\rm Im} \ll \, {\rm Ad}_a X_{\alpha}, \,
p_{\frak k} {\rm Ad}_a Y_{\alpha} \gg\\
& = & {\frac{2i}{\epsilon}} {\rm Im} \ll a^{\alpha} \ea - a^{-\alpha} \eb, \,
\, a^{-\alpha} (i \ea + i \eb) \gg \\
& = & {\frac{2i}{\epsilon}} (1\, - \, a^{-2 \alpha}).
\eeqa
Let $\sigma_x $ be the (left)-dressing vector
field on $AN$ defined by $x \in \fk$, i.e.,
$\sigma_x  = - x^l \backl \pi_{\scriptscriptstyle AN}$.
Then, taking $a = e^{-\lambda}$, we have
\beqa
\pi_{\scriptscriptstyle AN}(a) & = & \sum_{\alpha \in \Sigma_{+}}
{\frac{1}{\pi_{\scriptscriptstyle AN}(X_{\alpha}^{l}, Y_{\alpha}^{l})}}
\sigma_{X_{\alpha}}(a) \wedge \sigma_{Y_{\alpha}}(a).\\
& = &
- {\frac{i \epe }{2}}
\sum_{\alpha \in \Sigma_{+}}
{\frac{1}{1-e^{2\alpha(\lambda)}}}
\sigma_{X_{\alpha}}(a) \wedge \sigma_{Y_{\alpha}}(a)  \, \in \,
\wedge^2 T_a (K \cdot a).
\eeqa
Identify $K/T$ with $K \cdot a$ by $kT \mapsto k \cdot a$, we get
\[
\pi_{\lambda} (eT) \, = \, -{\frac{i \epe }{2}}
 \sum_{\alpha \in \Sigma_{+}}
{\frac{1}{1-e^{2\alpha(\lambda)}}}
 \Xa \wedge \Ya.
\]
Thus $\pi_{\lambda}$ is given as by (\ref{eq_pi-lambda}).
\qed

\subsection{$\pix$ as the result of Poisson induction}
\label{sec_induction}

We now look at the general case of $\pix$.
Set
\[
\fk_{\tx} \, = \, \ft \, + \,
{\rm span}_{\Bbb R} \{\Xa, \, \Ya: \, \a \in [X] \cap \Sigma_+\},
\]
and let $K_{\tx} \subset K$ be the connected subgroup of $K$ 
with Lie algebra $\fk_{\tx}$. We will
show that $\lix$ can be obtained
via Poisson induction  (see Remark \ref{rem_induction} below)
from a Poisson
structure on the smaller space $K_{\tx} /T$.

To this end, consider 
\[
\fk_{\tx}^{0} \, = \, \{\xi \in \fk^*: \,
\xi(x) = 0 \, \forall x \in \fk_{\tx} \}.
\]
Since we are identifying $\fk^*$ with $\fa + \fn$, we have
$\fk_{\tx}^{0}  \cong   \fn_{\tx}$ as real Lie
algebras, where $\fn_{\tx}$ is given in (\ref{eq_mx}).
Since $\fn_{\tx} \subset \fa + \fn$ is an ideal, we know
that $K_{\tx} \subset K$ is a Poisson subgroup \cite{lu-we:poi}.
In fact, set
\[
\Lambda_1 \, = \, -{\frac{i \epe }{2}} \sum_{\alpha \in [X] \cap \Sigma_{+}}
{\Xa \wedge \Ya \over 2}, \hspace{.3in} 
\Lambda_2 \,  = \,  -{\frac{i \epe }{2}} \sum_{\alpha \in \Sigma_{+}
\backslash [X]}
{\Xa \wedge \Ya \over 2}.
\]
Then, we have

\begin{prop}
\label{prop_poi-on-kx}
1) For any $x \in \fk_{\tx}, \, {\rm ad}_{x} \Lambda_2 \, = \, 0$;

2) The Poisson structure on $K_{\tx}$ (as a Poisson submanifold of 
$K$) is given by
\[
\pi_{\scriptscriptstyle K_X} (k_1) \, = \, R_{k_1} 
\Lambda_1 \, - \, L_{k_1} \Lambda_1,
\]
where $R_{k_1}$ and $L_{k_1}$ are respectively the right and left 
translations on $K_{\tx}$ by $k_1 \in K_{\tx}$.

3) The Manin triple for the Poisson Lie group 
$(K_{\tx}, \, \pi_{\scriptscriptstyle K_X})$ is 
$(\fm_{\tx}, \, \fk_{\tx}, \, \fa + \fu_{\tx}, \, 
{\frac{2i}{\epe}} \ll \, , \, \gg)$, where $\fm_{\tx}$,
given in (\ref{eq_gx}), is considered as over ${\Bbb R}$,
and $ \fu_{\tx} ={\rm span}_{\Bbb R} 
\{\ea, i \ea: \, \a \in [X] \cap \Sigma_+ \}$.
\end{prop} 

\noindent
{\bf Proof.} 1) Using the embedding of $\wedge^{\bullet} \fk $
into $\wedge^{\bullet} \fg$ as a real subspace, it is enough
to show that ${\rm ad}_{x} \Lambda_2 = 0$ for $x = 
E_{\alpha}$
with $\a \in [X]$.  Let $\alpha \in [X]\cap\Sigma_+$. Then,
\[
{\frac{2}{\epe}} {\rm ad}_{E_{\alpha}} \Lambda \, = \, 
\sum_{\beta \in \Sigma_{+}\backslash [X]} 
[\ea, \, E_{\beta}] \wedge E_{-\beta}
\, + \, E_{\beta} \wedge [\ea, \, E_{-\beta}].
\]
Set
\[
Y_{1}  =  \{\beta \in \Sigma_{+} \backslash [X]:  \a + \beta \in \Sigma\}, 
\hspace{.2in} {\rm and} \hspace{.2in} 
Y_2  =  \{\beta \in \Sigma_+ \backslash [X]:  \beta - \a \in \Sigma \}.
\]
Since $Y = \Sigma_{+} \backslash [X]$ has the property that
if $\a \in [X] \cap \Sigma_+$ and $\beta \in Y$ are such that
$\a + \beta \in \Sigma$ then $\a + \beta \in Y$, 
the map $ Y_1 \rightarrow  Y_2:  \beta \mapsto \a + \beta$
is a bijection. Thus
\beqa
{\frac{2}{\epe}} {\rm ad}_{E_{\alpha}} \Lambda_2
 & = & \sum_{\beta \in {\scriptscriptstyle Y_1}}( 
[\ea, E_{\beta}] \wedge  E_{-\beta}  +  
 E_{\alpha + \beta} \wedge [\ea,  E_{-(
\alpha + \beta)}])\\
 &=&  \sum_{\beta \in {\scriptscriptstyle Y_1}} 
(N_{\alpha, \beta} + N_{\alpha, -(\alpha + \beta)})
E_{\alpha + \beta} \wedge E_{-\beta} \\
 &=&  0.
\eeqa
Similarly, ${\rm ad}_{E_{-\alpha}} \Lambda_2 = 0$. 
This proves 1).

2) By definition, the induced Poisson structure 
$\pi_{\scriptscriptstyle K_X}$ on $K_{\tx}$ is the restriction of 
$\piK$
to $K_{\tx}$. Using the definition of $\piK$ and 1), we know
that $\pi_{\scriptscriptstyle K_X}$ is as given.

3) From the general theory of Poisson Lie groups \cite{lu-we:poi},
we know that the induced Lie algebra structure on $\fk_{\tx}^*$
is isomorphic to the quotient Lie algebra $\fk^* / \fk_{\tx}^{0}$.
Through the identifications $\fk^* \cong \fa + \fn$ and
$\fk_{\tx}^{0} \cong \fn_{\tx}$ via 
${\frac{2i}{\epe}} \ll \, , \, \gg$, we get 
$\fk_{\tx}^* \cong \fa + \fu_{\tx}$ via ${\frac{2i}{\epe}} \ll \, , \, \gg$
which is now considered as a symmetric scalar product on $\fm_{\tx}$
by restriction.
\qed

\begin{nota}
\label{nota_on-KXT}
{\em
Let $X_1 \subset X$ and let $\lambda  = \lambda_1 + 
{\pi i \over 2} \check{\rho}_{\txo} \in \fa_{\tx} +
{\pi i \over 2} \check{\rho}_{\txo}$ be such
that $\alpha(\lambda_1) \neq 0$ for
any $\alpha \in [X]$ with $\alpha(\check{\rho}_{\txo})$ even. 
By replacing $K$ by $K_{\tx}$
and by regarding $X$ as the set of all simple roots for
the root system for $(K_{\tx}, T)$, we know that
there is a $(K_{\tx}, \, \piKX)$-homogeneous
Poisson structure on  $K_{\tx}/T$ corresponding to
$X, X_1$ and $\lambda$. We will denote it by 
$\pk$.
}
\end{nota}

We now show that the Poisson structure $\pix$ on $K/T$ can
be obtained via Poisson induction from the Poisson structure
$\pk$ on $\kt / T$. 

To this end, consider the product space $K \times (\kt/T)$ with the 
product Poisson structure $\piK \oplus \pk$. Even though 
the diagonal (right) action of $\kt$ on $K \times (\kt / T)$
given by $k_1: (k, k^{'}T) \mapsto (kk_1, k_{1}^{-1} k^{'}T)$
 is in general
not Poisson, there is nevertheless a unique Poisson structure on
the quotient space $\kk$ such that the projection map
\[
K \times  (\kt / T) \lrw \kk: \, \, (k, \, k^{'}T) \Map
[(k, \, k^{'}T)]
\]
is a Poisson map. We temporarily denote this Poisson structure on
$\kk$ by $\pi_0$. 

\begin{rem}
\label{rem_induction}
{\em
In general, suppose that $K$ is a 
Poisson Lie group and $K_1 \subset K$ is a Poisson subgroup. Suppose that
$M$ is a Poisson manifold on which there is a Poisson action by $K_1$.
Then there is a unique Poisson structure on $K 
\times_{\scriptscriptstyle K_1} M$ such that the natural
projection from $K \times M$ to $K \times_{\scriptscriptstyle K_1} M$
is a Poisson map. Moreover, the left action of $K$ on 
$K \times_{\scriptscriptstyle K_1} M$ by left translations on the first 
factor is a Poisson action. We call this procedure of producing the Poisson
$K$-space $K \times_{\scriptscriptstyle K_1} M$ from the Poisson
$K_1$-space $M$ {\it Poisson induction}.
}
\end{rem}

\begin{prop}
\label{prop_induced}
We have $F_* \pi_0 = \pix$, where $F$ is the 
identification
\[
F: \, \, \kk \, \stackrel{\sim}{\lrw}\, K/T: \, \, 
[(k, \, k^{'}T)] \Map k k^{'}T.
\]
\end{prop}

\noindent
{\bf Proof.} Recall that $\pix$ is the image of $\tlp_{r_{\tx}(\lambda)}
= \Lambda^R - A_{\tx}(\lambda)^L$ 
under the projection $p_1: K \rightarrow K/T$, where $\Lambda^R$ (resp. 
$A_{\tx}(\lambda)^L$) is the right (resp. left) invariant bivector field
on $K$ with value $\Lambda$ (reps. $A_{\tx}(\lambda)$) at $e$,
and $A_{\tx}(\lambda) \in \fk \wedge \fk$ is the
skew symmetric part of
the $r$-matrix $r_{\tx}(\lambda)$ given in (\ref{eq_rx}).
On the other hand, $\pi_0$ is the image of $
\piK \oplus \bar{\pi}$ under the projection
\[
p_2: \, \, K \times \kt \lrw \kk: \, \, 
(k, \, k^{'}) \Map [(k, \, k^{'}T)],
\]
where $\bar{\pi}$ is the bi-vector field on $\kt$ defined by
$\bar{\pi}  = \Lambda_{1}^{R}  - \Lambda_{3}^{L}$ with 
\[
\Lambda_3 \, = \, -{\frac{i \epe }{2}} \sum_{\alpha \in [X] \cap \Sigma_+}
\coth \a (\lambda) {\Xa \wedge \Ya \over 2}.
\]
Because of the commutative diagram:
\beqa
K \times \kt \, \, \, \, \, \, \, \, \, \, \, \, &
\stackrel{m}{\lrw} &\, \, \, \, \, \, \, \, \, \, \, \,  K \\
\, \, \, \, \, \, \, \, \, \, \,  \, \, \,  \, \,  \, \, p_2 \downarrow 
 \, \, \, \, \, \, \, \, \, 
\, \, \, \, \, \, \, \, \, \, \, \, \, \, & & \, \, \, \, \, \, 
 \, \, \, \, \, \, \,
\downarrow p_1 \\
\, \kk \, \, \, \, \, \, & \stackrel{\lrw}{\scriptstyle F}&
  \, \, \, \, \, \, \, \, \, K/T,
\eeqa
where $ m: K \times \kt \lrw K: (k, \, k^{'}) \mapsto kk^{'}$,
it is enough to show that 
$m_{*} (\piK \oplus \bar{\pi})  =  
\tlp_{r_{\tx}(\lambda)}$,
or
\[
\tlp_{r_{\tx}(\lambda)} (kk_1) \, = \, L_k \bar{\pi}(k_1) \, + \, 
R_{k_1} \piK(k), \, \, \forall k \in K, \, k_1 \in \kt.
\]
But this follows easily from the definitions and the fact that
${\rm Ad}_{k_1} \Lambda_2 = \Lambda_2$ for all $k_1 \in \kt$.
\qed

We state some more properties of $\pix$ which can be proved either by
definitions or as corollaries of Proposition \ref{prop_induced}.

\begin{prop}
\label{prop_more-on-pix}
1) The embedding  $(\kt /T, \, \pk) \hookrightarrow (K/T, \, 
\pix)$ is a Poisson map;

2) With the Poisson structure $\piK$ on $K$, the 
Poisson structure $\pk$ on $\kt /T$ and the Poisson
structure $\pix$ on $K/T$, the map
\[
m_1: \, K \times (\kt /T) \lrw K/T: \, \, (k, \, k^{'}T) \Map kk^{'}T
\]
is a Poisson map;

3) Let $p_{*} \piK$ be the projection to $K/K_{\tx}$ of $\piK$
by $p: K \rightarrow K/K_{\tx}: k \mapsto kK_{\tx}$.
Then
the projection map $(K/T, \, \pix) \rightarrow (K/K_{\tx}, \, 
p_{*} \piK)$ is a Poisson map.
\end{prop}

\begin{rem}
\label{rem_bruhat}
{\em
The Poisson structure $p_{*} \piK$ on $K/K_{\tx}$ is known as the 
Bruhat-Poisson structure, because its symplectic leaves are
exactly the Bruhat cells in $K/K_{\tx}$. See \cite{lu-we:poi}.
}
\end{rem}

\subsection{The symplectic leaves of $\pix$}
\label{sec_leaves-1}

In this section, we  first describe the symplectic leaves of
$\pix$ for any $X \subset S(\Sigma_+)$ but  $X_1 = \emptyset$. 
The description of symplectic leaves for general
$\pix$ is somewhat complicated, and we will leave it
to the future. However, we will show that each $\pix$, for any
$X, X_1$ and $\lambda$, has at least one open symplectic
leaf.

\begin{nota}
\label{nota_pix-empty}
{\em
We will use $\qix$ to denote the Poisson structure
$\pix$ when $X_1$ is the empty set.
}
\end{nota}

We first recall that the space $K/T$ has the well-known Bruhat 
decomposition: Because of the  Iwasawa decomposition $G = KAN$
of $G$, the natural map $K/T \rightarrow G/B: kT \mapsto kB$
is a diffeomorphism. Its inverse map is $G/B \rightarrow K/T: 
gB \mapsto kT$ if $g = kan$ is the Iwasawa decomposition of $g$.
Thus we have
\[
K/T \, \cong \, G/B \, = \, \bigcup_{w \in W} N w B
\]
as a disjoint union. The set $N w B$ is called
the {\it Bruhat (or Schubert) cell} 
corresponding to $w \in W$. We denote it by $\sw$.
For $w \in W$, set
\[
\Phi_w \, = \, (-w\Sigma_{+}) \cap \Sigma_{+} \, = \, 
\{\a \in \Sigma_{+}: \, w^{-1} \a \in - \Sigma_{+} \}.
\]
Set
$\fn_w  =  {\rm span}_{\Bbb C} \{ \ea: \, \a \in \Phi_w \}$
and $ N_w  =  \exp \fn_w.$
Then $\sw$ is parametrized by $N_w$ by the map
\[
j_w: \, N_w \lrw \sw: \, \, n \Map n w B.
\]
Define
\beqa
j_1 & = & G \lrw K: \, \, g = kb \Map k \hspace{.2in}
{\rm for} \, \, k \in K, \, b \in AN;\\
j_2 & = & G \lrw K: \, \, g = bk \Map k \hspace{.2in}
{\rm for} \, \, k \in K, \, b \in AN.
\eeqa
Then we have a left action of $G$ on $K$ by
\[
G \times K \lrw K: \, \, (g, k) \Map g \circ k :=
j_1(gk),
\]
and a right action of $G$ on $K$:
\[
K \times G \lrw K: \, \, (k, g) \Map k^g := j_2(kg).
\]
The parametrization of $\sw$ by $N_w$ is then also given by
\[
j_w: \, N_w \lrw \sw: \, \, n \Map (n \circ \dw)T,
\]
where $\dw \in K$ is any representative of $w$ in $K$.

\begin{nota}
\label{nota_k-G}
{\em
For $k \in K$ and a subgroup $G_1 \subset G$, we set
\[
G_1 \circ k \, = \, \{g \circ k: \, g \in G_1 \} , 
\hspace{.5in}
k^{G_1} \, = \, \{k^g: \, g \in G_1 \}.
\]
}
\end{nota}

It is easy to show that 
$(AN) \circ k  =  k^{AN}$ for any $k \in K$. This set
is the symplectic leaf of $\piK$ in $K$ through the 
point $k$ (see \cite{soi:compact} \cite{lu-we:poi}).
Since $\kt \subset K$ is a Poisson submanifold,
we know that $(AN) \circ k  =  k^{AN} \subset \kt$ for $k \in \kt$.
 Moreover, 
if $w \in W$ and if
$\dw \in K$ is a representative of $w$ in $K$, set
\[
\cw \, = \, (AN) \circ \dw \subset K. 
\]
Then 
\begin{equation}
\label{eq_set-same}
\cw \, = \, (AN) \circ \dw \, =  \,N \circ \dw \, = \, N_w \circ \dw \, = \,
\dw^{AN} \,   =  \, \dw^N \, = \,  \dw^{N_{w^{-1}}}.
\end{equation}
Its image under the projection $K \rightarrow K/T$ is
the Bruhat cell $\sw$, which is also the symplectic leaf
of the Bruhat Poisson structure $\pi_{\infty}$ in $K/T$. See
\cite{soi:compact} \cite{lu-we:poi}.

Let $X \subset S(\Sigma_{+})$.
Denote by $\wx$ the subgroup of $W$ generated by the simple reflections
corresponding to elements in $X$. It is the 
Weyl group for $(\fm_{\tx}, \fh)$.
Introduce the subset $\wm$ of $W$:
\[
W^{\tx} \, = \, \{ w \in W: \, \, \Phi_{w^{-1}} \subset \Sigma_{+}\backslash
[X] \}.
\]
It follows from the definition that
$w \in \wm$ if and only if 
$w([X] \cap \Sigma_+) \subset \Sigma_{+}$. Moreover, we have $C_{\dot{w}_1}
= \dw_{1}^{N_{\tx}}$ for $w_1 \in \wm$ because $N_{w_{1}^{-1}} \subset
N_{\tx}$, where $N_{\tx}= \exp \fn_{\tx}$ with $\fn_{\tx}$ given by
(\ref{eq_mx}).
The following Lemma says that 
 each $w_1 \in \wm$ is the
minimal length representative for the coset $w_1 \wx$,
and that the set $\wm$ is a ``cross section" for the
canonical projection from $W$ to the coset space $W/\wx$.
For a proof of the Lemma, see \cite{ko:63}, Prop. 5.13.

\begin{lem}
\label{lem_minimal-rep}
For any $w \in W$, there exists a unique $w_1 \in \wm$ and
$w_2 \in \wx$ such that $w = w_1 w_2$.
Moreover, 
\[
\Phi_{w^{-1}} \, = \, \Phi_{w_{2}^{-1}} \cup w_{2}^{-1} \Phi_{w_{1}^{-1}}
\]
is a disjoint union, and the components on the right hand side are the
respective intersections of $\Phi_{w^{-1}}$ with $[X]$ and $\Sigma_{+}
\backslash [X]$.
Hence, $l(w) = l(w_1) + l(w_2)$.
\end{lem}

We can now describe the symplectic leaves of $\qix$ in $K/T$.

\begin{thm}
\label{thm_leaves-of-qix}
1) For each $w_1 \in \wm$, the union
$\bigcup_{w_2 \in W_{\tx}} \swb$ is the 
symplectic leaf of $\qix$ in $K/T$
through the point $w_1 \in K/T$.

2) These are all the symplectic leaves of $\qix$ in $K/T$.
\end{thm}

\noindent
{\bf Proof.} Set
\[
L_{\tx, \lambda} \, = \, e^{\lambda} \kt \el N_{\tx} \, = \,
N_{\tx} e^{\lambda} \kt \el.
\]
It is the connected  subgroup of $G$ with Lie algebra
\[
\fl_{\tx, \lambda} \, = \, {\rm Ad}_{e^{\lambda}}
\left(\fn_{\tx} \, + \, \fk_{\tx}\right).
\]
Notice that each $l \in 
L_{\tx, \lambda}$ can be written as a unique product $l = n_{\tx}
e^{\lambda} k \el$ for $n_{\tx} \in N_{\tx}$ and $k \in \kt$.

Denote by $S_{w_1}$ the symplectic leaf of $\qix$ through the point
$w_1 \in K/T$. Pick a representative
$\dw_1 $ of $w_1$ in $K$. 
By Theorem 7.2 of \cite{lu:homog}  (see also
\cite{ka:leaves}), 
the symplectic leaf $S_{w_1}$ is the image of the  set
$\dw_{1}^{L_{\tx, \lambda}}$ under the projection 
$K \rightarrow K/T$. We define a map
\[
M: \, \, L_{\tx, \lambda} \lrw N_{w_{1}^{-1}} \times \kt
\]
as follows: 
For $l = n_{\tx} e^{\lambda} k \el \in L_{\tx, \lambda}$, write
$k \el = b k^{'}$, where $b \in AU_{\tx}$ 
with $U_{\tx} = \exp \fu_{\tx}$ and $k^{'} \in \kt$, so that
$l = n_{\tx} e^{\lambda} b k^{'}$. Since 
the map $N_{w_{1}^{-1}} \rightarrow \cwa: n \mapsto
\dw_{1}^{n}$ is a diffeomorphism,
there exists a unique $n^{'} \in N_{w_{1}^{-1}}$ such that
$\dw_{1}^{n^{'}} = \dw_{1}^{n_{\tx} e^{\lambda} b}$. Now define 
$M(l) = (n^{'}, \, k^{'})$. It is easy to see that the map $M$ is onto
and that $\dw_{1}^{l} = \dw_{1}^{n^{'}} k^{'} \in \cwa \kt$.
This shows that
\[
\dw_{1}^{L_{\tx, \lambda}} \, = \, \cwa \kt.
\]
It is easy to show that the map
\[
\cwa \times \kt \lrw \cwa \kt: \, \, (c, \, k) \Map ck
\]
is a diffeomorphism, and that the image of $\cwa \kt$
to $K/T$ under the projection $K \rightarrow K/T$
is the union $\bigcup_{w_2 \in W_{\tx}} \Sigma_{w_1 w_2}$,
which is thus the symplectic leaf of the 
Poisson structure $\qix$ through the point $w_1 \in K/T$.
Now since   
\[
K/T \, = \, \bigcup_{w_1 \in W^{\tx}} S_{w_1}
\]
is already a disjoint union, we conclude that the collection 
$\{ S_{w_1}: \, w_1 \in W^{\tx} \}$ is that of all
symplectic leaves of $\qix$ in $K/T$.
\qed

Let  $w_1 \in W^{\tx}$. 
The following proposition identifies
the symplectic manifold  $S_{w_1}  
= \bigcup_{w_2 \in W_{\tx}} \swb$, as
a symplectic leaf of $\qix$ in $K/T$, with the product of two
symplectic manifolds. 
Recall that  for $w \in W$ with a representative $\dw$ in $K$,
the set $\cw \subset K$ is the symplectic leaf
of $\piK$ through the point $\dw$. Recall also from Notation
\ref{nota_on-KXT} the definition of the Poisson structure
$\pi_{\emptyset, \lambda}^{\tx}$ on $\kt/T$.
Note that it is symplectic by Proposition \ref{prop_dressing}.

\begin{prop}
\label{prop_product}
Let $w_1 \in \wm$ and  let $\dw_1$ be a representative
of $w_1$ in $K$. Equip $\cwa$ with the symplectic structure 
as a symplectic leaf of $\piK$ in $K$; Equip $\kt/T$ with
the symplectic structure $\pi^{\tx}_{\emptyset, \lambda}$,
and finally, 
equip $S_{w_1}$ with the 
symplectic structure as a symplectic leaf of 
$\qix$. Then the map
\[
m_1: \, \, \cwa \times \kt /T \lrw S_{w_1}: \, \, 
(k, \, k^{'}T) \Map k k^{'}T
\]
is a diffeomorphism between symplectic manifolds.
\end{prop}

\noindent
{\bf Proof.} This is a direct consequence of 
2) in Proposition \ref{prop_more-on-pix}.
\qed

Among all the elements in $\wm$, there is one which is the longest.
We denote this element by $w^{\tx}$, so $l(w^{\tx}) \geq l(w_1)$
for all $w_1 \in \wm$.

\begin{prop}
\label{prop_open-dense}
The symplectic leaf $S_{w^{\tx}}$ of $\qix$ in $K/T$ through the point
$w^{\tx}$ is
open and dense.
\end{prop}

\noindent
{\bf Proof.} Consider the projection $ K/T \rightarrow
K/\kt: kT \mapsto k\kt$. The image of $\Sigma_{w^{\tx}} \subset
K/T$ under this projection is an open dense subset (in fact a cell) in
$K/\kt$. Since $K/T \rightarrow K/\kt$ is a fibration, we know that
$S_{w^{\tx}}$ is open and dense in $K/T$.
\qed

\begin{cor}
\label{cor_finite}
Each Poisson structure $\qix$ has a finite number of symplectic leaves with
at least one of them open and dense.
\end{cor}

\begin{rem}
\label{rem_not-true}
{\em
Note that the statement in Corollary \ref{cor_finite} may not 
be true if $X_1 \neq \emptyset$, as is seen from case 3 of Example
\ref{exam_sl2}.
}
\end{rem}

The description of the symplectic leaves of $\pix$ in general 
is somewhat complicated. However, we have

\begin{prop}
\label{prop_non-degenerate}
The Poisson structure $\pix$ for $X = S(\Sigma_{+})$
(and $X_1 \subset X$ arbitrary) 
is non-degenerate at every element in the Weyl group $W$ of $(K, T)$
considered as a point in $K/T$. Consequently, the
symplectic leaves of $\pix$ through these points are open.
\end{prop}

\noindent
{\bf Proof.} Let $w \in W$ and let $\dot{w} \in K$ be a
representative of $w$ in $K$. Recall from the definition of
$\pix$ that $\pix = p_* \tilde{\pi}_1$, where
$p: K \rightarrow K/T$ is the natural projection
and $\tilde{\pi}_1$ is the bi-vector field on $K$ defined by
\[
\tilde{\pi}_1 \, = \, \Lambda^R \,  - \, A^L,
\]
with $\Lambda = -{i \epe \over 4} \sum_{\alpha \in \Sigma_{+}}
\Xa \wedge \Ya$ and
\[
A \, = \, -{i \epe \over 4} \sum_{\alpha \in \Sigma_{+}}
{e^{2\alpha(\lambda)} + 1 \over e^{2\alpha(\lambda)} -1}
\Xa \wedge \Ya.
\]
Thus
\beqa
l_{\dot{w}^{-1}} \tilde{\pi}_1(\dot{w}) & = & {\rm Ad}_{\dot{w}^{-1}}
\Lambda \, - \, A \\
& = & -{i \epe \over 4} \left(\sum_{\alpha \in \Sigma_{+}}
(X_{w^{-1}\alpha} \wedge Y_{w^{-1}\alpha}) \, + \,
\sum_{\alpha \in \Sigma_{+} }
({e^{2\alpha(\lambda)} + 1 \over e^{2\alpha(\lambda)} -1}
\Xa \wedge \Ya) \right)\\
& = & -{i \epe \over 4} \sum_{\alpha \in \Sigma_{+}, w \alpha < 0}
(1 + {e^{2\alpha(\lambda)} + 1 \over e^{2\alpha(\lambda)} -1})
\Xa \wedge \Ya \\
& & \hspace{.1in}\, - \, {i \epe \over 4}
\sum_{\alpha \in \Sigma_{+}, w \alpha > 0}
(-1 + {e^{2\alpha(\lambda)} + 1 \over e^{2\alpha(\lambda)} -1})
\Xa \wedge \Ya.
\eeqa
Since ${e^{2\alpha(\lambda)} + 1 \over e^{2\alpha(\lambda)} -1}
\neq \pm 1, \,
l_{\dot{w}^{-1}} \pix(\dot{w}T) = p_*
l_{\dot{w}^{-1}} \tilde{\pi}_1(\dot{w}) \in \wedge^2 T_e(K/T)$ is
non-degenerate. Hence
$\pix$ is non-degenerate at $w = \dot{w}T \in K/T$.
\qed

\begin{cor}
\label{cor_open-leaf}
For any $X, X_1$ and $\lambda$, the 
Poisson structure $\pix$ on $K/T$ 
has at least one open symplectic leaf.
\end{cor}

\noindent
{\bf Proof.} We use Proposition \ref{prop_induced} which says that
$\pix$ can be obtained via Poisson induction from the Poisson 
structure $\pk$ on $\kt/T$.  Recall the definition of
$\pk$ from Notation \ref{nota_on-KXT}. Since $X$ is the set of
all simple roots for the root systems for $(\kt, T)$, we know 
from Proposition \ref{prop_non-degenerate} that 
$\pk$ is non-degenerate at every Weyl group element in $W_{\tx}$,
regarded as points in $\kt/T$. Let $w_2 \in W_{\tx}$. Recall
that $w^{\tx}$ is the longest element in the set $W^{\tx}$. Let
$\dw^{\tx}$ be any representative of $w^{\tx}$ in $K$. Recall
that $C_{\dot{w}^{\tx}}$ is the symplectic leaf of $\piK$ in $K$
through $\dw^{\tx}$. By Proposition  \ref{prop_more-on-pix},
the map 
\[
(C_{\dot{w}^{\tx}}, \, \piK) \times (\kt/T, \pk) \, \lrw \, 
(K/T, \pix): \, 
(k, k^{'}T) \Map kk^{'}T
\]
is a Poisson map. But this map is a diffeomorphism onto its image which
is open because it is the inverse image under the natural
projection $K/T \rightarrow \kt/T$ of the biggest cell in $\kt/T$.
Thus the symplectic leaf of $\pix$ through the point 
$\dw^{\tx} w_2  \in K/T$ is open.
\qed

Note that the proof of Corollary \ref{cor_open-leaf} shows that
$\pix$ is open at every point in the coset $w^{\tx} W_{\tx} \subset K/T$.

\begin{exam}
\label{exam_sl2-again}
{\em  Corollary \ref{cor_open-leaf} can be
checked directly for the case of $\fg = \fsl(2, \C)$
by looking at the explicit formulas in Example \ref{exam_sl2}.
}
\end{exam}

\subsection{The modular vector fields and the leaf-wise moment maps for
the $T$-actions}
\label{sec_modular}

For an orientable Poisson manifold $(P, \pi)$ and a given volume form
$\mu$ on $P$,  the modular vector field of $\pi$ 
associated to $\mu$ is defined 
to be the vector field $v_{\mu}$ on $P$ satisfying
 $v_{\mu} \backl \mu = d(\pi \backl \mu)$. It measures
how Hamiltonian flows on $P$ fail to  preserve $\mu$. More details
can be found in \cite{we:modular}.

Coming back to $(K, \piK)$-homogeneous Poisson structures on $K/T$, we set
$\rho = {\frac{1}{2}} \sum_{\alpha \in \Sigma_{+}}
\alpha$ for the choice of $\Sigma_{+}$ in the definition of $\piK$.
Then we have $i H_{\rho} \in \ft$. We use
$\sigma_{i H_{\rho}} $ to denote the infinitesimal
generator of the $T$ action on $K/T$ by left translations in
the direction of $iH_{\rho}$.
 
\begin{prop}
\label{prop_modular}
For the Poisson structure $\piK$ on $K$ defined by (\ref{eq_on-K}) with
$\Lambda$ given in (\ref{eq_lambda-u}), all $(K, \piK)$-homogeneous
Poisson structures on $K/T$, and in particular all the $\pix$'s,
 have the same modular vector field $v$,
namely $v = -i \epe \sigma_{i H_\rho}$, with respect to a
(and thus any) $K$-invariant volume form on $K/T$.
\end{prop}

\begin{rem}
\label{rem_most-gneral}
{\em
Proposition \ref{prop_modular} is a statement about any Poisson Lie group
structure on $K$ since the Poisson structure $\piK$ on $K$
defined by (\ref{eq_on-K}) with
$\Lambda$ given in (\ref{eq_lambda-u}) is the most general form of
such structures.
}
\end{rem}
 
\noindent
{\bf Proof of Proposition \ref{prop_modular}.}
Let $\pi$ be an arbitrary
$(K, \piK)$-homogeneous Poisson structure. Then we know that
$\pi$ is the sum
\[
\pi \, = \, \pi(e)^{L} \, + \, p_* \piK,
\]
where $\pi(e)^{L}$ is the $K$-invariant bi-vector field on $K/T$ whose
value at $e = eT$ is $\pi(e)$, and $p_* \piK$ is
the projection of $\piK$ from $K$ to $K/T$ by $p: K \rightarrow
K/T: k \mapsto kT$ (it is the Bruhat Poisson
structure $\pi_{\infty}$ when $u = 0$ in the definition of $\Lambda$).
Let $\mu$ be a $K$-invariant volume form on $K/T$.
Let $b_{\mu}$ be the degree $-1$ operator on $\chi^{\bullet}(K/T)$
defined by $b_{\mu}(U) = (-1)^{|U|} d(U \backl \mu)$, 
so that $v = b_{\mu}(\pi)$ \cite{e-l-w:modular}. Then
$b_{\mu}(\pi) = b_{\mu}(\pi(e)^{L}) +
b_{\mu}(p_* \piK)$. Since $\mu$ is $K$-invariant, the operator
$b_{\mu}$ maps a $K$-invariant multi-vector field to another such.
Hence $b_{\mu}(\pi(e)^{L}) $ must be a $K$-invariant ($1$-)vector field
so it must be zero. Thus $b_{\mu}(\pi) = b_{\mu}(p_* \piK)$.
It is proved in \cite{e-l-w:modular}
that $b_{\mu}(p_* \piK) =
-i\epe \sigma_{i H_{\rho}}$, which is therefore the
modular vector field for any $\pi$.
\qed
 
The modular vector field is always a Poisson
vector field \cite{we:modular}, but it is not necessarily Hamiltonian
in general. For the rest of this section, we study this problem
for the modular vector field $v = -i \epe \sigma_{iH_{\rho}}$
for the Poisson structure $\qix$. We will show
that although $v$ is not globally Hamiltonian
unless $X = S(\Sigma_{+})$, it is leaf-wise, and we describe
its Hamiltonian function on each leaf. In fact, since
every $\qix$ is $T$-invariant (for the $T$-action on $K/T$
by left translations), we will 
describe the moment map for the $T$-action on each symplectic leaf
of $\qix$.
We are particularly
interested in the behavior of these moment maps when $\lambda$
goes infinity in various directions as in Section \ref{sec_limits}.
 
We first look at the Bruhat Poisson structure $\pi_{\infty}$
corresponding to $X = \emptyset$. This case (when $\epe = i$)
is studied in detail in \cite{lu:coor}. We recall the results there.
Denote by
\[
P_{\ta}: \, G = KAN \lrw A:\, \, g = kan \Map a,
\]
where $G = KAN$ is the Iwasawa decomposition of $G$ (as a real Lie
group).  For each $w \in W$, choose a representative
$\dw \in K$ of $w$ in $K$, and use
\[
j_w: \, \, N_w \lrw \sw: \, \, n \Map (n \circ \dw)T
\]
to parametrize the Bruhat cell $\sw$. For $n \in N_w$, let
$a_w(n) = P_{\ta}(n \dw) \in A$.  The element
$a_w(n)$ is independent of the choice of $\dw$, so we
have a well-defined map
\[
a_w: \, \, N_w \lrw A: \, \, n \Map a_w(n).
\]
Denote by $\Omega_w$ the symplectic structure on $\sw$ as a symplectic
leaf of $\pinf$. Then each $(\sw, \, \Omega_w)$ is a
Hamiltonian $T$-space.  The following fact is proved in \cite{lu:coor}.
 
\begin{prop}
\label{prop_bruhat-mom}
The map
\[
\phi_w: \, \sw \lrw \ft^*: \, \,
\la \phi_w, \, x \ra (kT) \, = \,
{2i \over \epe} {\rm Im} \ll {\rm Ad}_{\dot{w}}
\log a_w(j_{w}^{-1}(kT)), \, \, x \gg, \hspace{.3in}
x \in \ft
\]
is the moment map for the $T$-action on $(\sw, \, \Omega_w)$
such that $\phi_w(w) = 0$.
\end{prop}
 
In \cite{lu:coor}, we have written down an explicit formula 
for $\phi_w$ in certain Bott-Samelson type coordinates
$\{z_1, \bar{z}_1, z_2, \bar{z}_2, ..., z_{l(w)}, \bar{z}_{l(w)} \}$.
It takes the form
\[
\la \phi_w, \, x \ra \, = \,
-{1 \over \epe}\sum_{j=1}^{l(w)} {2\alpha_j (x) \over \ll \alpha_j,
\alpha_j \gg} \log(1+|z_j|^2) 
\]
where $\{\alpha_1, \alpha_2, ..., \alpha_{l(w)} \} = \Sigma_{+}
\cap (-w \Sigma_{+})$.
In particular, let $x = -i \epe (iH_{\rho}) =
\epe H_{\rho}$, we get a Hamiltonian function for
the vector field $v= -i \epe \sigma_{iH_{\rho}}$
on $(\sw, \Omega_w)$ as
\[
\la \phi_w, \, \epe H_\rho \ra \, = \,
-\sum_{j=1}^{l(w)} {2 \ll \rho, \alpha_j\gg \over \ll \alpha_j,
\alpha_j \gg} \log(1+|z_j|^2).
\]
 This function goes to $-\infty$ as $|z_j| \rightarrow \infty$
which corresponds to the boundary of $\sw$. Thus, the modular vector
field $v$ can not be globally Hamiltonian on $K/T$.

\bigskip
Next, we look at the case when $X = S(\Sigma_{+})$, so 
$\pi_{\tx, \emptyset, \lambda} = \pi_\lambda$ is the
the symplectic structure
on $K/T$ obtained by identifying
$K/T$ with the dressing orbit in the group $AN$
through the point $e^{-\lambda}$ (see Proposition \ref{prop_dressing}).
Since $K/T$ is simply connected,
The $T$-action on $K/T$ is Hamiltonian.
The following fact is proved in \cite{lu-ra:convexity}.
 
\begin{prop}
\label{prop_hamil}
The moment map for the $T$-action on $(K/T, \pi_{\lambda})$
is given by
\[
\Phi_{\lambda}: \, K/T \lrw \ft^*: \, \,
\la \Phi_{\lambda}, \, x \ra (kT) \, = \,
{2i \over \epe} {\rm Im} \ll \log (P_{\ta}(
k \el k^{-1})), \, \, x \gg, \hspace{.3in} x \in \ft.
\]
\end{prop}
 
\begin{rem}
\label{rem_convex}
{\em This fact plays the key role in the symplectic
proof of Kostant's nonlinear convexity theorem
given in \cite{lu-ra:convexity}.
}
\end{rem}
 
Corresponding to the fact that
$\lim_{t \rightarrow +\infty} \pi_{\lambda + t\check{\rho}} =
\pinf$, where $\check{\rho}$ is the sum of all
fundamental  coweights,
the two moment maps are related as follows.

\begin{prop}
\label{prop_limit-mom}
For any $\lambda \in \fa, w \in W$ and $kT \in \sw$,
\beqa
& & \lim_{t \rightarrow +\infty} \left( 
\Phi_{\lambda + t\check{\rho}} (kT) - \Phi_{\lambda + t\check{\rho}} 
(w) \right) 
\, = \, \phi_w(kT) \\
& & \lim_{t \rightarrow +\infty} d \Phi_{\lambda + t\check{\rho}} (kT)
\, = \, d \phi_w(kT).
\eeqa
\end{prop}

\noindent
{\bf Proof.} Using the parametrization of $\sw$ by
$N_w$, we regard both $\Phi_{\lambda + t\check{\rho}} |_{\Sigma_w}$
and $\phi_w$ as ($\ft^*$-valued) functions on $N_w$. Let
$n \in N_w$ with $k = n \circ \dw$. Write
\[
n \dw \, = \, k a_w(n) m(n)
\]
with $m(n) \in N_w$. Then
\[
e^{-\lambda} k^{-1} \, = \, 
(e^{-\lambda} a_w(n) m(n) a_w(n)^{-1} e^{\lambda} \dw^{-1} )
(\dw e^{-\lambda} a_w(n) \dw^{-1}) n^{-1}.
\]
Thus, for any $x \in \ft$,
\beqa
& & \la \Phi_{\lambda + t\check{\rho}} (n) - \Phi_{\lambda + t
\check{\rho}} (e) - 
\phi_w(n), \, \,  x\ra  \\
& = & 
{2i \over \epe} {\rm Im} \ll \log P_A (
e^{-\lambda -t\check{\rho}} a_w(n) m(n) a_w(n)^{-1} 
e^{\lambda + t\check{\rho}} \dw^{-1} ), \, \, x \gg,
\eeqa
where $e \in N_w$ is the identity element. Consider now the map
\[
\psi_t: \, N_w \lrw N_w: \, m \Map 
e^{-\lambda -t\check{\rho}} m e^{\lambda + t\check{\rho}}.
\]
Under the identification of $\fn_w$ with $N_w$ by the exponential map 
of $N_w$, this is the linear map ${\rm Ad}_{-\lambda -t\check{\rho}}$
on $\fn_w$, which goes to $0$ as $t \rightarrow +\infty$. Thus
\[
\lim_{t \rightarrow +\infty} \psi_t (m) = 0, \hspace{.2in}
{\rm and} \hspace{.2in}
\lim_{t \rightarrow +\infty} d\psi_t (m) = 0
\]
for all $m \in N_w$. But we have the composition of maps
\[
\la \Phi_{\lambda + t\check{\rho}} (n)  - \Phi_{\lambda + t
\check{\rho}} (e) -  
\phi_w(n), x \ra  \, = \, \eta_x (\psi_t(\xi(n))),
\]
where $\eta_x: N_w \rightarrow {\Bbb R}: m \mapsto {2i \over \epe}
{\rm Im} \ll \log P_A(m \dw^{-1}), x \gg$ and 
$\xi: N_w \rightarrow N_w: n \mapsto a_w(n) m(n) a_w(n)^{-1}$.
Thus the two limits in Proposition \ref{prop_limit-mom}
hold.
\qed

Now consider the general case of $\qix$.
Recall that the symplectic leaves of $\qix$ in $K/T$ are
indexed by elements in $W^{\tx}$. We keep the notation in
Proposition \ref{prop_product}, in which we have used the map $m_1$
to identify the symplectic leaf $S_{w_1}$ of $\qix$ in $K/T$ with
the product symplectic manifold $C_{\dot{w}_1} \times
K_{\tx}/T$. We use the projection map $C_{\dot{w}_1} \rightarrow 
\Sigma_{w_1}:
k \mapsto kT$ to identify $C_{\dot{w}_1}$ and $\Sigma_{w_1}$. This
identification is $T$-equivariant if we equip $C_{\dot{w}_1}$
with the $T$-action
\[
T \times C_{\dot{w}_1} \lrw C_{\dot{w}_1}: \, \, 
t \cdot k \Map t k (\dot{w}_{1}^{-1} t^{-1} \dw_1).
\]
Equip $C_{\dot{w}_1} \times K_{\tx}/T$ with the $T$-action
\[
T \times (C_{\dot{w}_1} \times K_{\tx}/T) \lrw
C_{\dot{w}_1} \times K_{\tx}/T: \, \, t \cdot (k, \, k^{'}T)
\Map (t k (\dw_{1}^{-1} t^{-1} \dw_1), \, \, \dw_{1}^{-1} t \dw_1 k^{'}T).
\]
Then the map $m_1$ in Proposition \ref{prop_product} is
$T$-equivariant. Denote by $\Phi_{\lambda, \tx}$
the moment map for the $T$-action on $(K_{\tx}/T, 
\pi_{\emptyset, \lambda}^{\tx})$. Then the moment map for the
$T$-action on $S_{w_1} \cong C_{\dot{w}_1} \times K_{\tx}/T$ is given by
\[
\la \phi_{\lambda, \tx, w_1} (k, \, k^{'}T), \, \, x \ra 
\, = \, \la \phi_{w_1}(kT), \, \, x \ra \, + \, 
\la \Phi_{\lambda, \tx}(k^{'}T), \, \, {\rm Ad}_{\dot{w}_{1}^{-1}} x \ra
\]
for all $x \in \ft$.

\begin{rem}
\label{rem_future}
{\em
There remain many problems to be addressed 
concerning the Poisson structures $\pix$. 
Other than the description of their 
symplectic leaves in the general case, one can try to compute its
Poisson cohomology according to the 
theory 
developed in \cite{lu:homog}. One can also study the 
$K$-invariant Poisson harmonic forms \cite{e-l:harm} of $\pix$. Another
problem is to construct the symplectic groupoids
for $\pix$. We hope to treat these problems in 
the future. 
}
\end{rem}

\end{document}